\documentclass[aap,MSNbibl,seceqn,citesort,dvips]{arximspdf}
\usepackage{mathbh}
\usepackage{dcolumn}
\usepackage{graphicx}

\doi{10.1214/11-AAP829} 
\volume{23}
\issue{2}
\pubyear{2013}
\firstpage{455}
\lastpage{491}

\makeatletter

\newcolumntype{d}[1]{D{.}{.}{#1}}

\newcommand{\one}{{\mathbf 1}}

\newcommand{\eqref}[1]{(\ref{#1})}

\newcommand{\cal}{\mathcal}

\newcommand{\esp}{\mathbb{E}} 

\newcommand{\indic}{\mathbh{1}} 

\newcommand{\proba}{\mathbb{P}}

\newcommand{\R}{\mathbb{R}}

\newcommand{\cF}{\mathcal{F}}

\newcommand{\bF}{\mathbb{F}}

\newcommand{\biota}{\bolds{\iota}}
\newcommand{\btau}{{\bolds{\tau}}}

\newcommand{\bsig}{\bolds{\sigma}}
\newcommand{\bL}{{\mathbf L}}
\newcommand{\btheta}{\bolds{\theta}}
\newcommand{\bl}{\bolds{\ell}}

\newcommand{\esssup}{\mathop{\operatorname{ess}\sup}}

\newcommand{\argmin}{\mathop{\arg\min}}

\newcommand{\I}{\mathbb{I}}

\newcommand{\Ee}{\mathbb{E}}

\newcommand{\F}{\mathbb{F}}

\renewcommand{\P}{\mathbb{P}}

\newcommand{\D}{\mathbb{D}}

\newcommand{\G}{\mathbb{G}}

\newcommand{\Sb}{{{\cal S}}}
\newcommand{\Lb}{\mathbf{L}}

\newcommand{\Ac}{{\cal A}}
\newcommand{\Bc}{{\cal B}}

\newcommand{\Dc}{{\cal D}}
\newcommand{\Ec}{{\cal E}}
\newcommand{\Fc}{{\cal F}}
\newcommand{\Gc}{{\cal G}}

\newcommand{\Pc}{{\cal P}}

\newcommand{\Oc}{{\cal O}}

\newproclaim{Def}{Definition}[section]
\newproclaim{Rem}{Remark}[section]
\newproclaim{Exe}{Example}[section]

\newtheorem{Thm}{Theorem}[section]

\makeatother

\begin{document}
\begin{frontmatter}

\title{Optimal investment under multiple defaults risk:
A~BSDE-decomposition approach}
\runtitle{Optimal investment under multiple defaults risk}

\begin{aug}
\author[A]{\fnms{Ying} \snm{Jiao}\ead[label=e1]{jiao@math.univ-paris-diderot.fr}},
\author[B]{\fnms{Idris} \snm{Kharroubi}\ead[label=e2]{kharroubi@ceremade.dauphine.fr}}
\and
\author[C]{\fnms{Huy\^en} \snm{Pham}\corref{}\ead[label=e3]{pham@math.univ-paris-diderot.fr}}
\runauthor{Y. Jiao, I. Kharroubi and H. Pham}
\affiliation{University Paris Diderot, University Paris Dauphine,
and University Paris Diderot, CREST-ENSAE and
Institut Universitaire de France}
\address[A]{Y. Jiao\\
Laboratoire de Probabilit\'es\\
\quad et Mod\`eles Al\'eatoires (LPMA)\\
University Paris Diderot\\
Site Chevaleret, Case 7012\\
75205 Paris Cedex 13\\
France\\
\printead{e1}} 
\address[B]{I. Kharroubi\\
CEREMADE\\
University Paris Dauphine\\
Place du Mar\'{e}chal DeLattre\\
\quad de Tassigny\\
75775 Paris Cedex 16\\
France\\
\printead{e2}}
\address[C]{H. Pham\\
LPMA\\
University Paris Diderot\\
Site Chevaleret, Case 7012\\
75205 Paris Cedex 13\\
and\\
CREST-ENSAE\\
\quad and Institut Universitaire de France\\
France\\
\printead{e3}}
\end{aug}

\received{\smonth{2} \syear{2011}}
\revised{\smonth{10} \syear{2011}}

%
\begin{abstract}
We study an optimal investment problem under contagion risk in a
financial model subject to multiple jumps and defaults. The global
market information is formulated as a progressive enlargement of a
default-free Brownian filtration, and the dependence of default times
is modeled by a conditional density hypothesis. In this It\^o-jump
process model, we give a decomposition of the corresponding stochastic
control problem into stochastic control problems in the default-free
filtration, which are determined in a backward induction. The dynamic
programming method leads to a backward recursive system of quadratic
backward stochastic differential equations (BSDEs) in Brownian
filtration, and our main result proves, under fairly general
conditions, the existence and uniqueness of a solution to this system,
which characterizes explicitly the value function and optimal
strategies to the optimal investment problem. We illustrate our
solutions approach with some numerical tests emphasizing the impact of
default intensities, loss or gain at defaults and correlation between
assets. Beyond the financial problem, our decomposition approach
provides a new perspective for solving quadratic BSDEs with a finite
number of jumps.
\end{abstract}

%

\end{frontmatter}

\section{Introduction}\label{sec1}

In this paper, we address an investment problem in an assets portfolio
subject to defaults and contagion risk, which is a major issue for risk
management in financial crisis period. We consider multiple default
events corresponding, for example, to the defaults of multi credit
names or to counter party defaults and contagion effects, meaning that
defaults on some assets may induce loss or gain on the other assets.
One usually formulates the default-free assets price process as an It\^
o process governed by some Brownian motion $W$, and
jumps are introduced at random default times, associated to a marked
point process $\mu$. The optimal investment problem in this incomplete
market framework may be then studied by stochastic control and dynamic
programming methods in the global filtration $\G$, generated by $W$ and
$\mu$. This leads in principle to Hamilton--Jacobi--Bellman
integrodifferential equations in a Markovian framework, and, more
generally, to backward stochastic differential equations (BSDEs) with
jumps, and the derivation relies on a martingale representation under
$\G$, with respect to $W$ and $\mu$, which holds under intensity
hypothesis on the defaults, and the so-called immersion property [or
(H) hypothesis]. Such an approach was used in the recent papers
\cite{ABE,LQ} in the single default case, and in \cite
{jeamatngo10} for the multiple defaults case. For exponential utility
criterion, the solution to the optimal investment problem is then
characterized through a quadratic BSDE with jumps, whose existence is
proved under a boundedness condition on the portfolio constraint set.

We revisit and extend the optimal investment problem in this multiple
defaults context by using an approach initiated in~\cite{JP09} in the
single default time case, and further developed in~\cite{Pham} in the
multiple defaults with random marks case. By viewing the global
filtration $\G$ as a progressive enlargement of filtrations of the
default-free filtration $\F$ generated by the Brownian motion $W$, with
the default filtration generated by the random times and jumps, the
basic idea is to split the global optimal investment problem, into
sub-control problems in the reference filtration $\F$ and corresponding
to optimal investment problems in default-free markets between two
default times.
More precisely, we derive a backward recursive decomposition by
starting from the
optimal investment problem when all defaults occurred, and then going
back to the initial optimal investment problem before any default.
The main point is to connect this family of stochastic control problems
in the $\F$-filtration, and this is achieved by assuming the existence
of a conditional density on the default times given the default-free
information $\F$. Such a density hypothesis, which is standard in the
theory of enlargement of filtrations, was recently introduced in
\cite{ejj1,ejj2} for credit risk analysis, and may be seen as an
extension of the usual intensity hypothesis.

This $\F$-decomposition approach allows us furthermore to formulate an
optimal investment problem where the portfolio constraint set can be
updated after each default time, depending possibly on the past
defaults, which is financially relevant. This extends the global
approach formulation where the portfolio set has to be fixed at the
beginning. Next, for exponential utility function criterion, we apply
dynamic programming method to each optimal investment problems in the
$\F$-filtration. We then get rid of the jump terms arising in the
dynamic programming in the $\G$-filtration, and are led instead to a
backward recursive system of quadratic BSDEs in Brownian filtration
with a nonstandard exponential term. Our main result is to prove under
fairly general conditions (without assuming in particular a boundedness
condition on the portfolio constraint set) the existence and uniqueness
of a solution to this system of BSDEs. Existence is showed by
induction, based on
Kobylanski results~\cite{kob00} together with approximating sequences
for dealing with the exponential term and unbounded portfolio, suitable
uniform estimates and comparison results for getting the convergence.
Uniqueness is obtained by verification arguments for relating the
solution of these BSDEs to the value functions of the $\F$-control
problems, and uses BMO-martingale tools. Moreover, an interesting
feature of our decomposition is to provide a nice characterization of
the optimal trading strategy between two default times, and to
emphasize the impact of defaults and jumps in the portfolio investment.
We also illustrate
numerically these results in a simple two defaultable assets model,
where each asset is subject to its own default and also to its counterpart.
Finally, we mention that beyond the optimal investment problem, the $\F
$-decomposition approach provides a new perspective for solving
(quadratic) BSDEs with finite number of jumps, see the recent paper
\cite{khalim11}.

The outline of this paper is organized as follows. In Section~\ref{sec2}, we
present the multiple defaults model where the assets price process is
written as a change of regimes model with jumps related to the default
times and random marks. Section~\ref{sec3} formulates the optimal investment
problem, and gives the decomposition of the corresponding stochastic
control problem. Section~\ref{sec4} is devoted to the derivation by dynamic
programming method of the sub-control problems in terms of a recursive
system of BSDEs, and to the existence and characterization results of
this system for the optimal investment problem. Finally, we provide in
Section~\ref{sec5} some numerical experiments for illustrating our solutions
approach in a simple two-defaultable assets model.

\section{Multiple defaults model}\label{sec2}

\subsection{Market information setup}\label{sec2.1}

We fix a probability space $(\Omega,\Gc, \proba)$,\break equipped with a
reference filtration $\mathbb F=(\Fc_t)_{t\geq0}$ satisfying the
usual conditions, and representing the default-free information
on the market. Let $\btau=(\tau_1,\ldots,\tau_n)$ be a
vector of $n$ random times, representing multiple default times,
and $\bL= (L_1,\ldots,L_n)$ be a vector of $n$ marks associated to
default times, $L_i$ being an $\Gc$-measurable random variable
taking values in some Polish space $E \subset\R^p$, and
representing, for example, the loss given default at time $\tau_i$.
The global market information is given by the default-free
information together with the observation of the default times
and their associated marks when they occur. It is then formalized
by the progressive enlargement of filtration $\G= \F\vee\D^1\vee\cdots
\vee\D^n$, where $\D^k = (\Dc_t^k)_{t\geq
0}$, $\Dc_t^k = \tilde\Dc_{t^+}^k$, $\tilde\Dc_t^k = \sigma(L_k 1_{\tau
_k\leq s},s\leq t) \vee\sigma(1_{\tau_k\leq s}, s \leq t)$, $k =
1,\ldots,n$. In
other words, $\G= (\Gc_t)_{t\geq0}$ is the smallest
right-continuous filtration containing $\F$ such that for any $k =
1,\ldots,n$, $\tau_k$ is a $\G$-stopping time, and $L_k$ is
$\Gc_{\tau_k}$-measurable.

For simplicity of presentation, we shall assume in the rest of this
paper that the default times are ordered, that is, $\tau_1\leq\cdots\leq
\tau_n$,
and so valued in $\Delta_n$ on $\{\tau_n<\infty\}$ where
\[
\Delta_k:= \{ (\theta_1,\ldots,\theta_k) \in(\R_+)^k\dvtx \theta_1\leq
\cdots\leq\theta_k \}.
\]
On one hand, this means that we do not distinguish specific credit
names, and only observe the successive default times, which is relevant
in practice for classical portfolio derivatives, like basket default
swaps. On the other hand, we may notice that the general nonordered
multiple random times case
for $(\tau_1,\ldots,\tau_n)$ [together with marks $(L_1,\ldots,L_n)$]
can be derived from the successive random times case by considering
suitable auxiliary marks. Indeed, denote by $\hat\tau_1\leq\cdots\leq
\hat\tau_n$ the corresponding ordered times, and by $\iota_k$ the index
mark valued in $\{1,\ldots,n\}$ so that $\hat\tau_k = \tau_{\iota
_k}$ for $k = 1,\ldots,n$. Then it is clear that the progressive
enlargement of filtration of $\F$ with the successive random times
$(\hat\tau_1,\ldots,\hat\tau_n)$, together with the marks
$(\iota_1,L_{\iota_1}, \ldots,\iota_n,L_{\iota_n})$, leads to the
filtration $\G$.

We introduce some notation used throughout the paper.
For any $(\theta_1,\ldots,\break\theta_n) \in\Delta_n$, $(\ell_1,\ldots
,\ell_n) \in E^n$,
we denote by $\btheta= (\theta_1,\ldots,\theta_n)$, $\bl= (\ell
_1,\ldots,\ell_n)$ and
$\btheta_k = (\theta_1,\ldots,\theta_k)$, $\bl_k = (\ell
_1,\ldots,\ell_k)$, for $k = 0,\ldots,n$, with the
convention $\theta_0 = \ell_0 = \varnothing$. We also denote by
$\btau_{k} = (\tau_1,\ldots,\tau_k)$ and $\bL_k = (L_1,\ldots,L_k)$.
For $t \geq0$, the set $\Omega_t^k$ denotes the event
\[
\Omega^{k}_{t}:= \{ \tau_k \leq t < \tau_{k+1} \}
\]
(with $\Omega_t^0 = \{t<\tau_1\}$, $\Omega_t^n = \{\tau_n\leq
t\})$ and represents the scenario where $k$ defaults occur before time $t$.
We call $\Omega_t^k$ as the $k$-default scenario at time~$t$. We define
similarly $\Omega_{t^-}^k = \{ \tau_k < t \leq\tau_{k+1}
\}$. Notice that for fixed $t$, the family $(\Omega_t^k)_{k=0,\ldots
,n}$ [resp., $(\Omega_{t^-}^k)_{k=0,\ldots,n}$] forms a partition of
$\Omega$.
We denote by $\Pc(\F)$ the $\sigma$-algebra of $\F$-predictable
measurable subsets on $\R_+\times\Omega$, and by
$\Pc_\F(\Delta^k,E^k)$ the set of indexed $\F$-predictable processes
$Z^k(\cdot,\cdot)$, that is, s.t. the map $(t,\omega,\btheta_k,\bl_k)
\rightarrow Z^k_t(\omega,\btheta_k,\bl_k)$ is $\Pc(\F)\otimes\Bc(\Delta
_k)\otimes
\Bc(E^k)$-measurable. We also denote by
$\Oc_\F(\Delta^k,E^k)$ the set of indexed $\F$-adapted processes
$Z^k(\cdot,\cdot)$, that is, such that for all $t \geq0$, the map
$(\omega
,\btheta_k,\bl_k) \rightarrow Z^k_t(\omega,\btheta_k,\bl_k)$ is $\Fc
_t\otimes\Bc(\Delta
_k)\otimes\Bc(E^k)$-measurable.

We recall from~\cite{Pham}, Lemma 2.1, or~\cite{Jiao09}, Lemma 4.1, the
key decomposition of any $\mathbb G$-adapted (resp., $\G$-predictable)
process $Z = (Z_t)_{t\geq0}$ in the form
\[
Z_t = \sum_{k=0}^n \indic_{\Omega_t^k} Z_t^k(\btau_k,\bL_k)
\qquad\Biggl[\mbox{resp., } Z_t = \sum_{k=0}^n \indic_{\Omega_{t^-}^k}
Z_t^k(\btau_k,\bL_k)\Biggr],\qquad t \geq0,
\]
where $Z^k$ lies in $\mathcal O_{\bF}(\Delta_k,E^k)$ [resp., $\mathcal
P_{\bF}(\Delta_k,E^k)$].

As in~\cite{ejj2} and~\cite{Pham}, we now suppose the existence of a
conditional joint density for $(\btau,\bL)$ with respect to the
filtration $\bF$.\vspace*{10pt}

\textit{Density hypothesis.}
There exists $\alpha\in\Oc_\F(\Delta_n,E^n)$ such that
for any bounded Borel function $f$ on $\Delta^n\times E^n$, and $t \geq0$,
%
\begin{equation} \label{density}
\esp[f(\btau,\bL)|\Fc_t] = \int_{\Delta^n\times E^n} f(\btheta,\bl
)\alpha_t(\btheta,\bl)\,d\btheta\eta(d\bl) \qquad\mbox{a.s.},
\end{equation}
where $d\btheta= d\theta_1\cdots d\theta_n$ is the Lebesgue
measure on $\R^n$, and $\eta(d\bl)$ is a Borel measure on $E^n$ in the form
$\eta(d\ell) = \eta_1(d\ell_1)\prod_{k=1}^{n-1}\eta_{k+1}(\bl
_k,d\ell_{k+1})$, with $\eta_1$ a nonnegative Borel measure on $E$ and
$\eta_{k+1}(\bl_k,d\ell_{k+1})$ a nonnegative transition kernel on
$E^k\times E$.
%
%
\begin{Rem} \label{remdens}
From condition \eqref{density}, we see that $\btau$ admits a
conditional (w.r.t.~$\F$) density with respect to the Lebesgue measure
given by $\alpha^{\btau}(\btheta) = \int\alpha(\btheta,\bl) \eta
(d\bl)$. This implies, in particular, that the default times are
totally inaccessible with
respect to the default-free information, which is consistent with the
financial modeling that the default events should arrive by
surprise, and cannot be read or predicted from the reference market
observation. This joint density condition w.r.t. the
Lebesgue measure also implies that the default times cannot occur
simultaneously, that is, $\tau_i \neq\tau_j$, $i\neq j$, a.s.,
which is a standard hypothesis in the modeling of multiple defaults.
Moreover, by considering a conditional density, and thus
a time-dependence of the martingale density process $(\alpha_t(\btheta
,\bl))_{t\geq0}$, we embed the relevant case in
practice when the default times are not independent of the reference
market information $\F$.
Compared to the classical default intensity processes for successive
defaults in the top-down modeling approach, the conditional density
provides more and necessary information for analyzing the impact of
default events. Further detailed discussion and some explicit models
for density of ordered random times are given in~\cite{ejj2}.

On the other hand, condition \eqref{density} implies that the family of
marks $\bL$ admits a conditional (w.r.t. $\F$) density with respect to
the measure $\eta(d\ell)$ given by $\alpha^{\bL}(\bl) = \int\alpha
(\btheta,\bl) \,d\btheta$. This general density hypothesis
\eqref{density} embeds several models of interest in applications. In
the case where $\alpha$ is separable in the form $\alpha(\btheta,\bl) =
\alpha^{\btau}(\btheta)\alpha^{\bL}(\bl)$, this means that the
random times and marks are independent given $\Fc_t$. The particular
case of nonrandom constant mark $L_k = \ell_k$ is obtained by
taking Dirac measure $\eta_k = \delta_{\ell_k}$. The case of
i.i.d. marks $L_k$, $k = 0,\ldots,n$, is included by taking
$\alpha^{\bL}(\bl)$ separable in $\ell_k$, and $\eta$ as a product
measure. We can also recover a density modeling of ordered default
times (as in the top-down approach) from a density model of the
nonordered defaults (as in the bottom-up approach). Indeed, let
$\btau=(\tau_1,\ldots,\tau_n)$ be a family of nonordered default times
having a density $\alpha^{\btau}$, and denote by $\hat\btau= (\hat\tau
_{1},\ldots,\hat\tau_{n})$, ${\biota} = (\iota_1,\ldots,\iota_n)$ the
associated ranked default times and
index marks. By using statistics order, we then see that
$(\btau,\biota)$ satisfy the density hypothesis with
\[
\hat\alpha(\theta_1,\ldots,\theta_n,i_1,\ldots,i_n) = \sum_{\bsig\in
\Sigma_n} \alpha^{\btau}\bigl(\theta_{\sigma(1)},\ldots,\theta_{\sigma(n)}\bigr)
1_{\{(i_1,\ldots,i_n)=(\sigma(1),\ldots,\sigma(n))\}}
\]
for $(\theta_1,\ldots,\theta_n) \in\Delta_n$, $\bl= (i_1,\ldots,i_n)
\in E = \{1,\ldots,n\}$, where $\Sigma_n$
denotes the set of all permutations $\bsig= (\sigma(1),\ldots
,\sigma(n))$
of $E$, and with $\eta(d\bl) = \sum_{\sigma\in\Sigma_n} \delta
_{\bl=\bsig}$,
$\eta_{k+1}(\bl_k,d\ell) = \sum_{i\in E\setminus\{\ell_1,\ldots
,\ell_k\}} \delta_{\ell=i}$.
\end{Rem}

\subsection{Assets and credit derivatives model}\label{sec2.2} \label{secasset}

We consider a portfolio of $d$ assets with value process defined by a
$d$-dimensional $\mathbb G$-adapted
process $S$. This process has the following decomposed form:
%
\begin{equation} \label{decS}
S_t = \sum_{k=0}^n \indic_{\Omega_t^k} S_t^{k}(\btau_k,\bL_k),
\end{equation}
where $S^{k}(\btheta_k,\bl_k)$, $\btheta_k = (\theta_1,\ldots,\theta
_k) \in\Delta_k$,
$\bl_k = (\ell_1,\ldots,\ell_k) \in E^k$, is an indexed process
in $\mathcal O_{\mathbb F}(\Delta_k,E^k)$,
valued in $\R_+^d$, representing the assets value in the $k$-default
scenario, given the past default events
$\btau_k = \btheta_k$ and the marks at default $\bL_k = \bl_k$.
Notice that $S_t$ is equal to the value $S_t^k$ only on the set $\Omega
_t^k$, that is, only for $\tau_k\leq t< \tau_{k+1}$.
We suppose that the dynamics of the indexed process $S^k$ is given by
%
\begin{equation}\label{SI}
dS_t^k(\btheta_k,\bl_k) = S_t^k(\btheta_k,\bl_k)*\bigl(b^k_t(\btheta_k,\bl
_k)\,dt+\sigma_t^k(\btheta_k,\bl_k)\,dW_t\bigr),\qquad
t \geq\theta_k,\hspace*{-35pt}
\end{equation}
where $W$ is a $m$-dimensional $(\P,\F)$-Brownian motion, $m \geq d$,
$b^k$ and $\sigma^k$ are indexed processes in $\mathcal P_{\mathbb
F}(\Delta_k,E^k)$, valued, respectively,
in $\R^d$ and $\R^{d\times m}$. Here, for $x = (x_1,\ldots,x_d)' \in\R
^d$ and
$y = (y_1,\ldots, y_d)'$ in $\R^{d\times q}$, the expression $x*y$
denotes the vector $(x_1y_1,\ldots,x_dy_d)'$ in $\R^{d\times q}$. Model
\eqref{decS}--\eqref{SI} can be viewed as an assets model with change of
regimes after each default event, with coefficients $b^k$, $\sigma^k$
depending on the past default times and marks. We make the usual
no-arbitrage assumption that there exists an indexed risk premium process
$\lambda^k \in\mathcal P_{\mathbb F}(\Delta_k,E^k)$ s.t. for all
$(\btheta_k,\bl_k) \in\Delta_k\times E^k$.
%
\begin{equation} \label{riskpremium}
\sigma_t^k(\btheta_k,\bl_k) \lambda_t^k(\btheta_k,\bl_k) = b_t^k(\btheta
_k,\bl_k),\qquad t \geq0.
\end{equation}
Moreover, in this contagion risk model, each default time may induce a
jump in the assets portfolio.
This is formalized by considering a family of indexed processes $\gamma
^k$, $k = 0,\ldots,n-1$, in $\Pc_{\F}(\Delta^k,E^k,E)$, and valued
in $[-1,\infty)^d$. For $(\btheta_k,\bl_k) \in\Delta^k\times E^k$,
and $\ell_{k+1} \in E$, $\gamma_t^k(\btheta_k,\bl_k,\ell_{k+1})$
represents the relative vector jump size on the $d$ assets at time $t =
\theta_{k+1} \geq\theta_k$
with a mark $\ell_{k+1}$, given the past default events $(\btau_k,\bL
_k) = (\btheta_k,\bl_k)$.
In other words, we have
%
\begin{equation} \label{sautS}
S^{k+1}_{\theta_{k+1}}(\btheta_{k+1},\bl_{k+1}) =
S^k_{\theta_{k+1}^-}(\btheta_k,\bl_k)*\bigl( {\one}_d + \gamma^k_{\theta
_{k+1}}(\btheta_k,\bl_k,\ell_{k+1}) \bigr),
\end{equation}
where we denote ${\one}_d$ as the vector in $\R^d$ with all components
equal to $1$.
%
\begin{Rem} \label{rempremium}
In this defaults market model, some assets may not be traded
anymore after default times, which means that their relative jump size
is equal to~$-1$. For $k = 0,\ldots,n$, $(\btheta_k,\bl_k) \in\Delta
_k\times E^k$, denote by $d^k(\btheta_k,\bl_k)$ the number of
assets among the
$d$-assets which cannot be traded anymore after $k$ defaults,
so that we can assume w.l.o.g. $b^k(\btheta_k,\bl_k) = (\bar
b^k(\btheta_k,\bl_k) 0)$,
$\sigma^k(\btheta_k,\bl_k) = (\bar\sigma^k(\btheta_k,\bl_k) 0
)$, $\gamma^k(\btheta_k,\bl_k,\ell) = (\bar\gamma^k(\btheta_k,\bl_k,\ell
) 0)$,
where $\bar b^k(\btheta_k,\bl_k)$, $\bar\sigma^k(\btheta_k,\bl_k)$,
$\bar\gamma^k(\btheta_k,\bl_k,\ell)$ are $\F$-predictable processes valued,
respectively, in $\R^{\bar d^k(\btheta_k,\bl_k)}$, $\R^{\bar d^k(\btheta
_k,\bl_k)\times m}$, $\R^{\bar d^k(\btheta_k,\bl_k)}$ with
$\bar d^k(\btheta_k,\bl_k) = d-d^k(\btheta_k,\bl_k)$, the number of
remaining tradable assets. Either $\bar d^k(\btheta_k,\bl_k) = 0$,
and so
$\sigma^k(\btheta_k,\bl_k) = 0$, $b^k(\btheta_k,\bl_k) = 0$,
$\gamma^k(\btheta_k,\bl_k,\ell) = 0$,
in which case \eqref{riskpremium} is trivially satisfied, or
$\bar d^k(\btheta_k,\bl_k) \geq1$, and we shall assume the natural
condition that the volatility matrix $\bar\sigma^k(\btheta_k,\bl_k)$
is of full rank. We can then define the risk premium
\[
\lambda^k(\btheta_k,\bl_k) = \bar\sigma^k(\btheta_k,\bl_k)'(\bar\sigma
^k(\btheta_k,\bl_k)\bar\sigma^k(\btheta_k,\bl_k)')^{-1}
\bar b^k(\btheta_k,\bl_k),
\]
which satisfies \eqref{riskpremium}.
\end{Rem}
%
\begin{Rem} \label{remdynmu}
One can write the dynamics of the assets model \eqref{decS}--\eqref
{SI}--\eqref{sautS} as a jump-diffusion process under $\G$. Let us
define the $\G$-predict\-able processes
$(b_t)_{t\geq0}$ and $(\sigma_t)_{t\geq0}$ valued, respectively, in
$\R^d$ and $\R^{d\times m}$ by
%
\begin{equation} \label{defbsig}
b_t = \sum_{k=0}^n \indic_{\Omega_{t^-}^k} b_t^k(\btau_k,\bL_k),\qquad
\sigma_t = \sum_{k=0}^n \indic_{\Omega_{t^-}^k} \sigma_t^k(\btau
_k,\bL_k),
\end{equation}
and the indexed $\G$-predictable process $\gamma$, valued in $\R^d$,
and defined by
\[
\gamma_t(\ell) = \sum_{k=0}^{n-1} \indic_{\Omega_{t^-}^k} \gamma
_t^k(\btau_k,\bL_k,\ell).
\]
Let us introduce the random measure $\mu(dt,d\ell)$ associated to the
jump times and marks $(\tau_k,L_k)$, $k = 1,\ldots,n$, and given by
%
\begin{equation} \label{defmu}
\mu([0,t]\times B) = \sum_k 1_{\tau_k\leq t} 1_{L_k \in B},\qquad t
\geq0, B \in\Bc(E).
\end{equation}
Then, the dynamics of the assets value process $S$ is written under $\G
$ as
%
\begin{equation} \label{dynS}
dS_t = S_t*\biggl( b_t \,dt + \sigma_t \,dW_t + \int_E \gamma_t(\ell) \mu
(dt,d\ell) \biggr).
\end{equation}
Notice that in formulation \eqref{dynS}, the process $W$ is not in
general a Brownian motion under $(\P,\G)$, but a semimartingale under
the density hypothesis, which preserves
the semimartingale property [also called (H$'$) hypothesis in the
progressive enlargement of filtrations literature]. We also mention
that the random measure $\mu$ is not independent of $W$ under the
conditional density hypothesis. Thus, in general, we de not have a
martingale representation theorem under $(\P,\G)$ with respect to $W$
and $\mu$.
\end{Rem}

In this market, a credit derivative of maturity $T$ is modeled by a $\Gc
_T$-measurable random variable $H_T$, thus decomposed in the form
%
\begin{equation} \label{decH}
H_T = \sum_{k=0}^n 1_{\Omega_T^k} H_T^k(\btau_k,\bL_k),
\end{equation}
where $H_T^k(\cdot,\cdot)$ is $\Fc_T\otimes\Bc(\Delta_k)\otimes\Bc
(E^k)$-measurable, and represents the option payoff when $k$ defaults occured
before maturity $T$.

The above model setup is quite general, and allows us to consider a
large family of explicit examples.

\subsection{Examples}\label{sec2.3}

\begin{Exe}[(Exogenous counterparty default)]
We consider a highly risky underlying name (e.g., Lehman Brothers)
which may have an impact on many
other names once the default occurs. One should take into consideration
this counterparty risk for each asset in the
investment portfolio; however, the risky name itself is not
necessarily contained in the investment portfolio. A special case of
this example containing one asset (without marks) has been considered
in~\cite{JP09}; see also~\cite{LQ,ABE}.

There is one default time $\tau$ ($n = 1$), which may induce jumps
in the price process $S$ of the $d$-assets portfolio.
The corresponding mark is given by a random vector $L$ valued in $E
\subset[-1,\infty)^d$, representing the proportional jump size in the
$d$-assets price.

The assets price process is described by
\[
S_t = S_t^{0} \indic_{t<\tau} + S_t^{1}(\tau,L) \indic_{t\geq\tau},
\]
where $S^{0}$ is the price process before default, governed by
\[
dS_t^{0} = S_t^{0}*(b_t^{0} \,dt + \sigma_t^{0} \,dW_t)
\]
and the indexed process $S^{1}(\theta, \ell)$, $(\theta,\ell)\in\R
_+\times E$, representing the price process after default at time
$\theta$ and with mark $\ell$, is given by
\begin{eqnarray*}
dS_t^{1}(\theta,\ell)&=&S_t^{1}(\theta,\ell)*\bigl(b_t^{1}(\theta,\ell) \,dt +
\sigma_t^{1}(\theta,\ell) \,dW_t\bigr),\qquad t\geq\theta,\\
S_\theta^{1}(\theta,\ell) &=& S_\theta^{0}*({\one}_d + \ell).
\end{eqnarray*}
Here $W$ is an $m$-dimensional $(\P,\bF)$-Brownian motion, $m \geq d$,
$b^{0}$, $\sigma^{0}$ are
$\bF$-predictable bounded processes valued, respectively, in $\R^d$ and
$\R^{d\times m}$, and the indexed processes
$b^{1}$, $\sigma^{1}$ lie in $\Pc_{\F}(\R_+,E)$, and valued,
respectively, in $\R^d$ and $\R^{d\times m}$.
\end{Exe}
%
\begin{Exe}[(Assets portfolio with multilateral counterparty risks)]
\label{Exebonds}
The defaults family and the assets family coincide, each underlying
name subjected to the default risk of
itself and to the counterparty default risks of the other names of the
portfolio.\vadjust{\goodbreak} The assets family is represented by a
portfolio of defaultable bonds. Recall that a defaultable bond is a
credit derivative which insures 1 euro to its buyer if no default
occurs before the maturity; otherwise,
the buyer of the bond receives a recovery rate at the default time. The
recovery rate may be random, and so it is viewed in our model as a
random mark at the default time.

In this contagion risk model, the number of defaults times $n$ is
equal to the number $d$ of defaultable bonds. We denote by $P^i$ the
price process of the $i$th defaultable bond of maturity $T_i$, by
$\tau_i$ its default time and $L_i$ its (random) recovery rate valued
in $E = [0,1)$. The price process $P^i$ drops to $L_i$ at the
default time $\tau_i$, and remains constant afterward. Moreover, at the
default times $\tau_j$, $j \neq i$ (which are not necessarily
ordered) of the other defaultable bonds, the price process $P^i$ has a
jump, which may depend on $\tau_j$ and $L_j$. Actually, the jump size
of $P^i$ will typically depend on $L_j$ if the name $i$ is the debt
holder of name $j$. The assets portfolio price process $S = (P^1,\ldots
,P^n)$ has the decomposed form
%
\begin{equation} \label{decSi} P^i_t =
\sum_{k=0}^n \indic_{\hat\tau_k\leq t< \hat\tau_{k+1}}
P_t^{i,k}(\hat\btau_k,\biota_k,\hat\bL_k),\qquad t \geq0,
\end{equation}
where
$\hat\btau_k = (\hat\tau_1,\ldots,\hat\tau_k)$ denotes the $k$
first ordered times, $\biota_k = (\iota_1,\ldots,\iota_k)$ the
corresponding index marks, that is, $\hat\tau_k = \tau_{\iota_k}$,
and $\hat\bL_k = (L_{\iota_1},\ldots,L_{\iota_k})$. The index
$\F$-adapted process $P^{i,k}(\btheta_k,\biota_k,\bl_k)$, for
$(\btheta_k,\biota_k,\bl_k) \in\Delta_k\times\I^k\times E^k$,
represents the price process of the $i$th defaultable bond, given that
the $k$ names $(\iota_1,\ldots,\iota_k)$ defaulted at times
$\hat\btau_k = \btheta_k$ with the marks $\hat\bL_k = \bl_k$.
Here, we denoted by $\I_k = \{(\iota_1,\ldots,\iota_k) \in
\{1,\ldots,n\}\dvtx \iota_j \neq\iota_{j'} \mbox{ for } j \neq j' \}$.
When $i \in\{\iota_1,\ldots,\iota_k\}$, that is, $i = \iota_j$
for some $j = 1,\ldots,k$, then $P^{i,k}(\btheta_k,\biota_k,\bl_k) =
\ell_j$, and otherwise it evolves according to the dynamics
\begin{eqnarray*}
&&
dP_t^{i,k}(\btheta_k,\biota_k,\bl_k)\\
&&\qquad =
P_t^{i,k}(\btheta_k,\biota_k,\bl_k)\bigl(
b_t^{i,k}(\btheta_k,\biota_k,\bl_k) \,dt +
\sigma_t^{i,k}(\btheta_k,\biota_k,\bl_k) \,dW_t\bigr),\qquad t \geq
\theta_k.
\end{eqnarray*}
Here $W$ is an $m$-dimensional $(\P,\bF)$-Brownian
motion, $m \geq n$, and the indexed processes $b^{i,k}$,
$\sigma^{i,k}$ lie in $\Pc_{\F}(\Delta_k,\I^k,E^k)$, and are valued,
respectively, in $\R^n$ and are $\R^{1\times m}$. The jumps of the
$i$th defaultable bond are given by
\[
P^{i,k+1}_{\theta_{k+1}}(\btheta_{k+1},\biota_{k+1},\bl_{k+1}) =
P^{i,k}_{\theta_{k+1}^-}(\btheta_k,\biota_k,\bl_k)
\bigl( 1 + \gamma^{i,k}_{\theta_{k+1}}(\btheta_k,\biota_k,\bl_k,\iota
_{k+1},\ell_{k+1}) \bigr)
\]
for $\theta_{k+1} \geq\theta_k$, and $\iota_{k+1} \in\{
1,\ldots,n\}\setminus\{\iota_1,\ldots,\iota_k\}$, and we have
$\gamma^{i,k}_{\theta_{k+1}}(\btheta_k,\biota_k,\bl_k$, $\iota_{k+1},\ell
_{k+1}) = -1 + \ell_{k+1}/P^{i,k}_{\theta_{k+1}^-}(\btheta_k,\biota
_k,\bl_k)$,
meaning that
$P^{i,k+1}_{\theta_{k+1}}(\btheta_{k+1},\biota_{k+1},\break\bl_{k+1}) = \ell
_{k+1}$, when $\iota_{k+1} = i$.
This model is compatible with several ones in the literature (see,
e.g.,~\cite{BC2008,CJZ}), and we shall focus in the last
section on this example for numerical illustrations in the case $n = 2$.
\end{Exe}


\begin{Exe}[(Basket default swaps)]\label{Exebasket}
A $k$th-to-default swap is a credit derivative contract, which provides\vadjust{\goodbreak}
to its buyer the protection against the $k$th default of the underlying name.
The protection buyer pays a regular {continuous} premium $p$ until the
occurrence of the $k$th default time, or until the maturity $T$, if
there are less than $k$ defaults before maturity. In return, the
protection seller pays the loss $1-L_k$ where $L_k$ is the recovery
rate if $\tau_k$ is the $k$th default occurring before $T$, and zero
otherwise. By considering that the available information consists in
the ranked default times and the corresponding recovery rates, and
assuming zero interest rate,
the payoff of this contract can then be written in the form \eqref
{decH} with
\[
H_T^i(\btheta_i,\bl_i) = \cases{
{ - p \theta_{k} + (1-\ell_{k})}, &\quad if $i \geq k$, \cr
- pT, &\quad if $i < k$,}
\]
for $\btheta_i = (\theta_1,\ldots,\theta_i) \in\Delta_i$, $\bl
_i = (\ell_1,\ldots,\ell_i) \in E^i$.\vspace*{-2pt}
\end{Exe}

\section{The optimal investment problem}\label{sec3}\vspace*{-2pt}

\subsection{Trading strategies and wealth process}\label{sec3.1}

A trading strategy in the $d$-assets portfolio model described in
Section~\ref{secasset} is a $\G$-predictable process $\pi$, hence
decomposed in the form
%
\begin{equation} \label{decpi}
\pi_t = \sum_{k=0}^n 1_{\Omega_{t^-}^k} \pi_t^k(\btau_k,\bL_k),\qquad
t \geq0,
\end{equation}
where $\pi^k$ is an indexed process in $\mathcal P_{\mathbb F}(\Delta_k,E^k)$,
and $\pi^k(\btheta_k,\bl_k)$ is valued in $A^k$ closed set of $\R^d$
containing the zero element, and represents
the amount invested continuously in the $d$-assets in the $k$-default
scenario, given the past default events $\btau_k = \btheta_k$ and
the marks at default $\bL_k = \bl_k$, for $(\btheta_k,\bl_k) \in\Delta
_k\times E^k$. Notice that in this modeling,
we allow the space $A^k$ of strategies constraints to vary between
default times.
This means that the investor can update her portfolio constraint set
based on the observation of the past default events,
and this includes the typical case for defaultable bonds where the
assets cannot be traded anymore after their own defaults.
Notice that this framework is then more general than the standard
formulation of a stochastic control problem, where
the control set $A$ is invariant in time.
%
\begin{Rem}
It is possible to formulate a more general framework for the
modeling of portfolio constraints by considering that the set $A^k$ may
depend on the past defaults and marks. More precisely, by introducing
for any $k = 0,\ldots,n$,
a~closed set $\bar A^k \subset\R^d\times\Delta_k\times E^k$, s.t.
$(0,\btheta_k,\bl_k) \in\bar A^k$ for all $(\btheta_k,\bl_k) \in\Delta
_k\times E^k$, and denoting by $A^k(\btheta_k,\bl_k) = \{
\pi\in\R^d\dvtx (\pi,\btheta_k,\bl_k) \in\bar A^k\}$, the portfolio
constraint is defined by the condition that the process $\pi^k(\btheta
_k,\bl_k)$ should be valued in $A^k(\btheta_k,\bl_k)$. In the rest of
this paper, and for simplicity of notation, we shall focus on the case
where $A^k$ does not depend on the past defaults and marks, that is,
$\bar A^k = A^k\times\Delta_k\times E^k$.
\end{Rem}

In the sequel,\vspace*{1pt} we shall often identify the strategy $\pi$ with the family
$(\pi^k)_{k=0,\ldots,n}$ given in \eqref{decpi}, and we require\vadjust{\goodbreak} the
integrability conditions: for all
$\btheta_k \in\Delta_k$, $\bl_k \in E^k$,
%
\begin{eqnarray} \label{integpi}
&&\int_0^T |\pi_t^k(\btheta_k,\bl_k)'b_t^k(\btheta_k,\bl_k)| \,dt\nonumber\\
&&\quad{}+ \int_0^T |\pi_t^k(\btheta_k,\bl_k)'\sigma_t^k(\btheta_k,\bl
_k)|^2 \,dt\\
&&\qquad < \infty \qquad\mbox{a.s.},\nonumber
\end{eqnarray}
where $T < \infty$ is a fixed finite horizon time.
Given a trading strategy $\pi= (\pi^k)_{k=0,\ldots,n}$, the
corresponding wealth process is defined by
%
\begin{equation} \label{defXpi}
X_t = \sum_{k=0}^n 1_{\Omega_{t}^k} X_t^{k}(\btau_k,\bL_k),\qquad
0\leq t \leq T,
\end{equation}
where $X^{k}(\btheta_k,\bl_k)$, $\btheta_k \in\Delta_k$, $\bl_k \in
E^k$, is an indexed process in
$\mathcal O_{\mathbb F}(\Delta_k,E^k)$, representing the wealth
controlled by $\pi^k(\btheta_k,\bl_k)$ in the price process $S^k(\btheta
_k,\bl_k)$,
given the past default events $\btau_k = \btheta_k$ and the marks
at default
$\bL_k = \bl_k$. From the dynamics \eqref{SI}, and under \eqref
{integpi}, it is governed by
%
\begin{equation} \label{dynXI}\qquad
dX_t^{k}(\btheta_k,\bl_k) = \pi_t^k(\btheta_k,\bl_k)'\bigl( b_t^k(\btheta
_k,\bl_k) \,dt + \sigma^k(\btheta_k,\bl_k) \,dW_t \bigr),\qquad
t \geq\theta_k.
\end{equation}
Moreover, each default time induces a jump in the assets price process,
and then also on the wealth process. From \eqref{sautS}, it is given by
%
\[
X^{k+1}_{\theta_{k+1}}(\btheta_{k+1},\bl_{k+1}) = X^{k}_{\theta
_{k+1}^-}(\btheta_k,\bl_k)
+ \pi_{\theta_{k+1}}^k(\btheta_k,\bl_k)'\gamma_{\theta_{k+1}}^k(\btheta
_k,\bl_k,\ell_{k+1}).
\]
Notice that the dynamics of the wealth process can be written as a
jump-It\^o controlled process under $\G$ by means of the random measure
$\mu$ in \eqref{defmu},
%
\begin{equation} \label{dynXG}
dX_t = \pi_t'\biggl(b_t \,dt + \sigma_t \,dW_t + \int_E \gamma_t(\ell) \mu
(dt,d\ell) \biggr).
\end{equation}

\subsection{Value functions and $\F$-decomposition}\label{sec3.2}

Let $U$ be an exponential utility with risk aversion coefficient $p > 0$,
\[
U(x) = - \exp( - p x),\qquad x \in\R.
\]
We
consider an investor with preferences described by the utility function
$U$, who can trade in the $d$-assets portfolio following an
\textit{admissible} trading strategy $\pi\in\Ac_{\G}$ to be defined
below, associated with a wealth process $X = X^{x,\pi}$, as in
\eqref{defXpi} with initial capital $X_{0^-} = x$. Moreover, the
investor has to deliver at maturity $T$ an option of payoff $H_T$, a
bounded $\Gc_T$-measurable random variable, decomposed into the form
\eqref{decH}. The optimal investment problem is then defined by
%
\begin{equation}
\label{defV0} V^0(x) = \sup_{\pi\in\Ac_{\G}} \Ee[
U(X_T^{x,\pi}-H_T) ].
\end{equation}
Our main\vspace*{1pt} goal is to provide existence
and characterization results of the value function $V^0$, and
of the optimal trading\vadjust{\goodbreak} strategy $\hat\pi$ (which does not depend on the
initial wealth $x$ from the exponential form of $U$) in the general
assets framework described in the previous section. A first step is to
define in a suitable way the set of admissible trading strategies.
%
\begin{Def}[(Admissible trading strategies)] \label{defadmi}
For $k = 0,\ldots,n$, $\Ac^k_{\F}$ denotes the set of
indexed process $\pi^k$ in
$\mathcal P_{\mathbb F}(\Delta_k,E^k)$, valued in $A^k$ satisfying
\eqref{integpi}, and such that:
\begin{itemize}
\item the family $\{U(X^{k}_{\tau}(\btheta_{k},\bl_{k})), \tau\
\F\mbox{-stopping time valued in }[\theta_{k},T]\}$ is uniformly
integrable, that is, $U(X^k(\btheta_k,\bl_k))$ is of class
(D);\vspace*{1pt}
\item$\Ee[ \int_{\theta_k}^T \int_{E} (-U)(X^{k}_{s}(\btheta
_{k},\bl_{k})+ \pi^k_{s}(\btheta_{k},\bl_{k})'\gamma^k_{s}(\btheta
_{k},\bl_{k},\ell))\eta_{k+1}(\bl_k,d\ell)\,ds ] < \infty$,\break
when $k \leq n-1$,
\end{itemize}
for all $(\btheta_{k},\bl_{k})\in\Delta_k(T)\times E^k$, where we set
$\Delta_k(T) = \Delta_k\cap[0,T]^k$.
We then denote by $\Ac_{\G} = (\Ac^k_{\F})_{k=0,\ldots,n}$ the set of
admissible trading strategies $\pi= (\pi^k)_{k=0,\ldots,n}$.
\end{Def}

As mentioned above, the indexed control sets $A^k$ in which the trading
strategies take values may vary after each default time.
This nonstandard feature in control theory prevents a direct resolution
to \eqref{defV0} by dynamic programming or duality methods in the global
filtration $\G$, relying on the dynamics \eqref{dynXG} of the controlled
wealth process.
Following the approach in~\cite{Pham}, we then provide a decomposition
of the global optimization problem \eqref{defV0} in terms of a family of
optimization problems with respect to the default-free filtration $\bF
$. Under the density hypothesis~\eqref{density}, let us define a
family of auxiliary processes $\alpha^k \in\Oc_\F(\Delta_k,E^k)$,
$k = 0,\ldots,n$, which is related to the survival probability and
is defined by recursive induction from $\alpha^n = \alpha$,
%
\begin{equation}\label{AlphaAuxk}\quad
\alpha_t^k(\btheta_k,\bl_k) = \int_t^\infty\int_E \alpha_t^{k+1}(\btheta
_k,\theta_{k+1},\bl_k,\ell_{k+1}) \,d\theta_{k+1} \eta_{k+1}(\bl
_k,d\ell_{k+1})
\end{equation}
for $k = 0,\ldots,n-1$, so that
\[
\P[ \tau_{k+1} > t | \Fc_t] = \int_{\Delta_k\times E^k} \alpha
_t^k(\btheta_k,\bl_k) \,d\btheta_k \eta(d\bl_k),\qquad
\P[ \tau_{1} > t | \Fc_t] = \alpha_t^0,
\]
where $d\btheta_k = d\theta_1\cdots d\theta_k$, $\eta(d\bl_k) = \eta
_1(d\ell_1)\cdots\eta_{k}(\bl_{k-1},d\ell_{k})$.
Given $\pi^k \in\Ac_{\F}^k$, we denote by $X^{k,x}(\btheta_k,\bl
_k)$ the controlled process solution to \eqref{dynXI} and starting from
$x$ at $\theta_{k}$. For simplicity of notation, we omit the dependence
of $X^{k,x}$ in $\pi^k$.
The value function to the global $\G$-optimization problem \eqref{defV0}
is then given in a backward induction from the
$\F$-optimization problems:
%
\begin{eqnarray}\qquad
\label{Vn}
&&
V^{n}(x,\btheta,\bl) \nonumber\\[-8pt]\\[-8pt]
&&\qquad= \esssup_{\pi^{n}\in\Ac^{n}_{\F}}
\Ee[ U(X_T^{n,x} - H_T^{n}) \alpha_T(\btheta,\bl) | \Fc
_{\theta_{n}} ]  ,\nonumber\\
\label{VI}
&&
V^k(x,\btheta_k,\bl_k) \nonumber\\
&&\qquad= \esssup_{\pi^{k}\in\Ac^{k}_{\F}} \Ee\biggl[
U(X_T^{k,x} - H_T^{k} ) \alpha_T^k(\btheta_k,\bl_k)
\nonumber\\[-8pt]\\[-8pt]
&&\hspace*{43pt}\qquad\quad{}+ \int_{\theta_{k}}^T \int_E V^{k+1}\bigl( X_{\theta_{{k+1}}}^{k,x} +
\pi_{\theta_{{k+1}}}^k \gamma_{\theta_{{k+1}}}^k(\ell_{{k+1}}),
\btheta_{{k+1}},\bl_{{k+1}}\bigr)\nonumber\\
&&\hspace*{154pt}\qquad\quad{}\times \eta_{{k+1}}(\bl_k,d\ell_{{k+1}})
\,d\theta_{{k+1}} \Big| \Fc_{\theta_{k}} \biggr] \nonumber
\end{eqnarray}
for any $x \in\R$, $k =0,\ldots,n$, $(\btheta_k,\bl_k) \in\Delta
_k(T)\times E^k$. Here $X^{k,x}$ denotes wealth process in \eqref
{dynXI} controlled by $\pi^k$, and starting from $x$ at time $\theta
_k$. To alleviate notation, we omit, and often omit in the sequel, in
$X^{k,x}$, $H_T^k$, $\pi^k$, $\gamma^k$, the dependence on $(\btheta
_k,\bl_k)$, when there is no ambiguity.
Notice that $(\btheta_k,\bl_k)$ appears in \eqref{VI} as a
parameter index through $X^{k,x}$, $H_T^k$, $\pi^k$, $\gamma^k$ and
$\alpha^k$. On the other hand, $\btheta_k$ appears also
via $\theta_{k}$ as the initial time in \eqref{VI}. The interpretation
of relations \eqref{Vn}--\eqref{VI} is the following. $V^k$ represents the
value function of the optimal investment problem in the $k$-default
scenario, and equality \eqref{VI} may be understood as a dynamic
programming relation between two consecutive default times: on the
$k$-default scenario, with a wealth controlled process $X^k$, either
there are no other defaults before time~$T$ (which is measured by the
survival density $\alpha^k$), in which case, the investor receives the
terminal gain $U(X_T^k-H_T^k)$, or there is a default at time $\tau
_{k+1}$, which occurs between $\theta_k$ and $T$, inducing a jump on
$X^k$, and from which the maximal expected profit is $V^{k+1}$.
Moreover, if there exists, for all $k = 0,\ldots,n$, some $\hat\pi
^k \in\Ac_\F^k$ attaining the essential supremum in \eqref
{Vn}--\eqref{VI}, then the trading strategy $\hat\pi= (\hat\pi
^k)_{k=0,\ldots,n} \in\Ac_\G$, is optimal for the initial
investment problem~\eqref{defV0}.

\section{Backward recursive system of BSDEs}\label{sec4}

In this section, we exploit the specific form of the exponential
utility function $U(x)$ in order to characterize, by dynamic
programming methods, the solutions to the stochastic optimization
problems \eqref{Vn}--\eqref{VI} in terms of a recursive system of indexed
backward stochastic differential equations (BSDEs) with respect to the
filtration $\F$, assumed from now on to be generated by the
$m$-dimensional Brownian motion~$W$.

We use a verification approach in the following sense. We first derive
formally the system of BSDEs associated to the $\F$-stochastic control
problems. The main step is then to obtain existence of a solution to
these BSDEs, and prove
that this BSDEs-solution indeed provides the solution to our optimal
investment problem.

Let us consider the starting problem \eqref{Vn} of the backward
induction. For fixed $(\btheta,\bl) \in\Delta_n(T)\times E^n$,
problem \eqref{Vn} is a classical exponential utility maximization in
the market model $S^{n}(\btheta,\bl)$ starting from $\theta_{n}$,
and with random endowment $\tilde H_T^{n} = H_T^{n}+\frac{1}{p}\ln
\alpha_T$.
We recall briefly how to derive the corresponding BSDE.
For $t \in[\theta_n,T]$, $\nu^n \in\Ac_\F^n$, let us
introduce the following set of controls coinciding with $\nu$ until
time $t$:
\[
\Ac_\F^n(t,\nu^n) = \{ \pi^n \in\Ac_\F^n\dvtx \pi_{\cdot\wedge t}^n = \nu
_{\cdot\wedge t}^n \}
\]
and define the dynamic version of \eqref{Vn} by considering the
following family of $\F$-adapted processes:
%
\begin{equation} \label{defVndyn}
V_t^n(x,\btheta,\bl,\nu^n) =
\esssup_{\pi^n\in\Ac_\F^n(t,\nu^n)} \Ee[ U(X_T^{n,x} - \tilde
H_T^{n}) | \Fc_{t} ],\qquad t \geq\theta_n,
\end{equation}
so that $V_{\theta_n}^n(x,\btheta,\bl,\nu^n) = V^n(x,\btheta,\bl)$
for any $\nu^n \in\Ac_\F^n$. From the dynamic programming principle,
one should have the supermartingale property of $\{ V_t^n(x,\btheta,\bl
,\nu^n)$, $\theta_n\leq t\leq T\}$, for any $\nu^n \in\Ac_\F^n$, and if
an optimal control exists for \eqref{defVndyn}, we should have the
martingale property of
$\{ V_t^n(x,\btheta,\bl,\hat\pi^n), \theta_n\leq t\leq T\}$ for some
$\hat\pi^n \in\Ac_\F^n$. Moreover, from the exponential form of
the utility function $U$ and the additive form of the wealth process
$X^n$ in \eqref{dynXI}, the value function process $V^n$ should be in
the form
\[
V_t^n(x,\btheta,\bl,\nu^n) = U\bigl(X_t^{n,x} - Y_t^n(\btheta,\bl)
\bigr),\qquad \theta_n \leq t \leq T,
\]
for some indexed $\F$-adapted process $Y^n$ independent of $\nu^n$,
that we search in the form: $dY_t^n = -f_t^ndt + Z_t^n \,dW_t$.
Then, by using the above supermartingale and martingale property of the
dynamic programming principle, and since
$V_T^n(x,\btheta,\bl,\nu^n) = U(x-\tilde H_T^n)$ by \eqref
{defVndyn}, we see that $(Y^n,Z^n)$ should satisfy the following
indexed BSDE:
{\renewcommand{\theequation}{En}
\begin{eqnarray}\label{En}
Y_t^{n}(\btheta,\bl) &=&
H_T^{n}(\btheta,\bl) +\frac{1}{p} \ln\alpha_T(\btheta,\bl)
\nonumber\\[-8pt]\\[-8pt]
& &{} + \int_t^T f^{n}(r,Z_r^{n},\btheta,\bl) \,dr - \int_t^T
Z_r^{n}\,dW_r, \qquad\theta_{n}\leq t\leq T,\nonumber
\end{eqnarray}}

\noindent and the generator $f^{n}$ is the indexed process in $\Pc_\F(\R^m,\Delta
_n,E^n)$ defined by
%
\setcounter{equation}{1}
\begin{eqnarray}\label{deffn}
f^{n}(t,z,\btheta,\bl) &=& \inf_{\pi\in A^{n}} \biggl\{ \frac{p}{2}
| z - \sigma_t^{n}(\btheta,\bl)'\pi|^2 - b^{n}(\btheta,\bl)'\pi
\biggr\}  \nonumber\\
&=& - \lambda_t^{n}(\btheta,\bl)z - \frac{1}{2p} |\lambda
_t^{n}(\btheta,\bl)|^2\\
&&{}+ \frac{p}{2} \inf_{\pi\in A^{n}} \biggl| z + \frac{1}{p} \lambda
_t^{n}(\btheta,\bl) - \sigma_t^{n}(\btheta,\bl)' \pi\biggr|^2, \nonumber
\end{eqnarray}
where the second equality comes from \eqref{riskpremium}. This quadratic
BSDE is similar to the one considered in~\cite{REK} or~\cite{HIM} in a
default-free market.
Next, consider the problems \eqref{VI}, and define similarly\vadjust{\goodbreak} the dynamic
version by considering the value function process
%
\begin{eqnarray}\label{defVkdyn}
&&
V_t^k(x,\btheta_k,\bl_k,\nu^k) \nonumber\\[-2pt]
&&\qquad= \esssup_{\pi^{k}\in\Ac^{k}_{\F}(t,\nu^k)}
\Ee\biggl[ U\bigl(X_T^{k,x} - H_T^{k}(\btheta_k,\bl_k) \bigr) \alpha
_T^k(\btheta_k,\bl_k) \nonumber\\[-9pt]\\[-9pt]
&&\hspace*{59pt}\qquad\quad{}+ \int_{t}^T \int_E V^{k+1}_{\theta_{k+1}} \bigl( X_{\theta
_{{k+1}}}^{k,x} + \pi_{\theta_{{k+1}}}^k \gamma_{\theta
_{{k+1}}}^k(\ell_{{k+1}}),
\btheta_{{k+1}},\bl_{{k+1}}\bigr)\nonumber\\[-2pt]
&&\qquad\quad\hspace*{174pt}{}\times \eta_{{k+1}}(\bl_k,d\ell_{{k+1}})
\,d\theta_{{k+1}} \Big| \Fc_{t} \biggr] \nonumber
\end{eqnarray}
for $\theta_k\leq t\leq T$, where $\Ac_\F^k(t,\nu^k) = \{ \pi^k
\in\Ac_\F^k\dvtx \pi_{\cdot\wedge t}^k = \nu_{\cdot\wedge t}^k \}$, for
$\nu
^k \in\Ac_{\F}^k$, so that $V_{\theta_k}^k(x,\btheta_k,\bl_k,\nu^k) =
V^k(x,\btheta_k,\bl_k)$.
The dynamic programming principle for \eqref{defVkdyn} formally implies
that the process
\begin{eqnarray*}
&&
V_t^k(x,\btheta_k,\bl_k,\nu^k)\\[-2pt]
&&\qquad{} + \int_0^t \int_E V^{k+1}\bigl( X_{\theta
_{k+1}}^{k,x} + \nu_{\theta_{k+1}}^k \gamma_{\theta_{k+1}}^k(\ell_{k+1}),
\btheta_{k+1},\bl_{k+1}\bigr)\\[-2pt]
&&\qquad\quad\hspace*{28.4pt}{}\times \eta_{k+1}(\bl_k,d\ell_{k+1}) \,d\theta_{k+1}
\end{eqnarray*}
for $\theta_k\leq t\leq T$ is a $(\P,\F)$-supermartingale for any $\nu
^k \in\Ac_{\F}^k$, and is a martingale for $\hat\pi^k$ if it is
an optimal control for \eqref{defVkdyn}.
Again, from the exponential form of the utility function $U$, the
additive form of the wealth process $X^k$ in \eqref{dynXI}, and by
induction, we see that the value function process $V^k$ should be in
the form
\[
V_t^k(x,\btheta_k,\bl_k,\nu^k) = U\bigl(X_t^{k,x} - Y_{t}^{k}(\btheta
_k,\bl_k)\bigr),\qquad \theta_k \leq t \leq T,
\]
for some\vspace*{1pt} indexed $\F$-adapted process $Y^{k}$, independent of $\nu^k$,
that we search in the form
$dY_t^k = -f_t^kdt + Z_t^k \,dW_t$. By using the supermartingale and
martingale properties of the dynamic programming principle for $V^k$,
and since $V_T^k(x,\btheta_k,\bl_k) = U(x-\tilde H_T^k)$, with
$\tilde H_T^{k} = H_T^{k}+\frac{1}{p}\ln\alpha_T^k$, we see that
$(Y^k,Z^k)$ should satisfy the indexed BSDE,
{\renewcommand{\theequation}{Ek}
\begin{eqnarray}\label{EK}\qquad
Y_t^{k}(\btheta_k,\bl_k) &=&
H_T^{k}(\btheta_k,\bl_k) +\frac{1}{p}\ln\alpha_T^k(\btheta_k,\bl_k)
+ \int_t^T f^{k}(r,Y_r^k,Z_r^{k},\btheta_k,\bl_k)
\,dr\nonumber\\[-9pt]\\[-9pt]
&&{} - \int_t^T Z_r^{k}\,dW_r, \qquad
\theta_{k} \leq t\leq T,\nonumber
\end{eqnarray}}

\noindent with a generator $f^k$ defined by
%
\setcounter{equation}{3}
\begin{eqnarray}\label{deffk}
&&
f^{k}(t,y,z,\btheta_k,\bl_k)\nonumber\\[-2pt]
&&\qquad= \inf_{\pi\in A^{k}} \biggl\{ \frac
{p}{2} | z - \sigma_t^{k}(\btheta_k,\bl_k)'\pi|^2 -
b_t^{k}(\btheta_k,\bl_k)'\pi
\nonumber\\
&&\hspace*{25.4pt}\qquad\quad{} + \frac{1}{p} U(y) \int_E U\bigl(\pi\gamma_t^k(\btheta
_k,\bl_k,\ell)\nonumber\\
&&\qquad\quad\hspace*{97.8pt}{} - Y_t^{k+1}(\btheta_k,t,\bl_k,\ell)\bigr)
\eta_{k+1}(\bl_k,d\ell) \biggr\} \nonumber\\[-8pt]\\[-8pt]
&&\qquad= - \lambda_t^{k}(\btheta_k,\bl_k)  z - \frac{1}{2p} |\lambda
_t^{k}(\btheta_k,\bl_k)|^2 \nonumber\\
&&\qquad\quad{} + \inf_{\pi\in A^{k}} \biggl\{ \frac{p}{2} \biggl| z + \frac
{1}{p} \lambda_t^{k}(\btheta_k,\bl_k) - \sigma_t^{k}(\btheta_k,\bl
_k)'\pi\biggr|^2 \nonumber\\
&&\qquad\quad\hspace*{38.2pt}{} + \frac{1}{p} U(y) \int_E U\bigl(\pi\gamma
_t^k(\btheta_k,\bl_k,\ell)\nonumber\\
&&\qquad\quad\hspace*{110.5pt}{}
- Y_t^{k+1}(\btheta_k,t,\bl_k,\ell)\bigr) \eta_{k+1}(\bl_k,d\ell) \biggr\}
, \nonumber
\end{eqnarray}
where the second equality comes from \eqref{riskpremium}.

The equations (\ref{EK}), $k = 0,\ldots,n$, define thus a recursive
system of families of BSDEs, indexed by $(\btheta,\bl) \in\Delta
_n(T)\times E^n$, and the rest of this section is devoted first
to the well-posedness and existence of a solution to this system, and
then to its uniqueness via a verification theorem relating the solution
to the value functions \eqref{defVndyn}, \eqref{defVkdyn}.

\subsection{Existence to the recursive system of indexed BSDEs}\label{sec4.1}

The generators of our system of BSDEs do not satisfy the usual
Lipschitz or quadratic growth assumptions.
In particular, in addition to the growth condition in $z$ for $f^k$
defined in \eqref{deffk},
there is an exponential term in $y$ via the utility function $U(y)$,
which prevents a direct application of known existence results in the
literature for BSDEs.

Let us introduce some notation for sets of processes. We denote by $\Sb
_c^\infty[t,T]$ the set of $\F$-adapted continuous processes $Y$ which
are essentially bounded on $[t,T]$, that is, $\|Y\|_{{\Sb_c^\infty
[t,T]}}:= \esssup_{(s,\omega)\in[t,T]\times\Omega}|Y_s(\omega)| <
\infty$,
and by $\Lb^2_W[t,T]$ the set of $\F$-predictable processes $Z$ s.t.
$\Ee[\int_t^T |Z_s|^2 \,ds] < \infty$.
For any $k = 0,\ldots,n$, we denote by $\Sb_c^\infty(\Delta_k,E^k)$
the set of indexed
$\F$-adapted continuous processes $Y^k$ in $\Oc_{\F}(\Delta_k,E^k)$,
which are essentially bounded, uniformly in their indices
\[
\|Y^k\|_{{\Sb_c^\infty(\Delta_k,E^k)}}:= \sup_{(\btheta_k,\bl_k)\in
\Delta_k(T)\times E^k}
\|Y^k(\btheta_k,\bl_k)\|_{{\Sb_c^\infty[\theta_k,T]}}
< \infty.
\]
%
We also denote by $\Lb^2_W(\Delta_k,E^k)$ the set of indexed $\F
$-predictable processes $Z^k$ in $\Pc_\F(\Delta_k,E^k)$ such that
\[
\Ee\biggl[ \int_{\theta_k}^T | Z^k_t(\btheta_k,\bl_k)|^2 \,dt \biggr] <
\infty\qquad\forall(\btheta_k,\bl_k) \in\Delta_k(T)\times E^k.
\]

We make the following boundedness assumptions:

{\renewcommand\thelonglist{(HB)}
\renewcommand\labellonglist{\thelonglist}
\begin{longlist}
\item\label{HB}
%
(i) The risk premium is bounded uniformly w.r.t. its indices: there
exists a constant $C > 0$
such that for any $k = 0,\ldots,n$, $(\btheta_k,\bl_k) \in\Delta
_k(T)\times E^k$, $t \in[\theta_k,T]$,
\[
|\lambda_t^k(\btheta_k,\bl_k)|
\leq C \qquad\mbox{a.s.}
\]

(ii) The indexed $\Fc_T$-measurable random variables $H_T^k$ and $\ln
\alpha_{T}^k$ are bound\-ed uniformly in their indices:
there exists a constant $C > 0$ such that for any $k = 0,\ldots
,n$, $(\btheta_k,\bl_k) \in\Delta_k(T)\times E^k$,
\[
|H_T^k(\btheta_k,\bl_k)| + |{\ln\alpha_T^k}(\btheta_k,\bl_k)| \leq C
\qquad\mbox{a.s.}
\]
\end{longlist}}

We then state the existence result for the recursive system of BSDEs.
%
\begin{Thm}\label{exsredsr}
Under~\ref{HB}, there exists a solution $(Y^k,Z^k)_{k=0,\ldots,n} \in
\prod_{k=0}^n \Sb_c^\infty(\Delta_k,E^k) \times\Lb^2_W(\Delta
_k,E^k)$ to the recursive system of indexed BSDEs (\ref{EK}), $k =
0,\ldots,n$.
\end{Thm}
\begin{pf}
We prove the result by a backward induction on $k = 0,\ldots,n$,
and consider the property
\renewcommand{\theequation}{$\Pc_k$}
\begin{equation}\label{Pck}
\mbox{there exists a solution $Y^k \in\Sb
_c^\infty(\Delta_k,E^k)$ to (\ref{EK}).}
\end{equation}

$\bullet$ For $k = n$. From expression \eqref{deffn} of the
generatof $f^{n}$, there exists some positive constant $C$ s.t.
\begin{eqnarray}
|f^{n}(t,z,\btheta,\bl)| \leq C\bigl(|z|^2+ |\lambda_{t}^{n}(\btheta,\bl
)|^2\bigr)\nonumber\\
&&\eqntext{\forall(t,z,\btheta,\bl)\in[0,T]\times\R^m\times\Delta_n(T)\times
E^{n}.}
\end{eqnarray}
Hence, under~\ref{HB}, we can apply Theorem 2.3 in~\cite{kob00}
for any fixed $(\btheta,\bl) \in\Delta_n(T)\times E^n$, and get
the existence of a solution $(Y^{n}(\btheta,\bl), Z^{n}(\btheta,\bl))
\in\Sb_c^\infty[\theta_n$, $T]\times\Lb_W^2[\theta_{n},T]$. Moreover,
from Proposition 2.1 in~\cite{kob00}, we have the following estimate:
\begin{eqnarray*}
|Y^{n}_{t}(\btheta,\bl)| &\leq&\esssup_{\Omega}\biggl(|H_{T}(\btheta,\bl
)|+{1\over p} |{\ln\alpha_{T}}(\btheta,\bl)|\biggr)\\
&&{}+ C \int_{t}^T|\lambda^{n}_{s}(\btheta,\bl)|^2\,ds,\qquad \theta_{n}
\leq t\leq T.
\end{eqnarray*}
Under~\ref{HB}, this implies that $\sup_{(\btheta,\bl)\in\Delta
_n(T)\times E^n} \|Y^{n}(\btheta,\bl)\|_{{\Sb_c^\infty[\theta_n,T]}} <
\infty$.
Finally, the measurability of $Y^n$ and $Z^n$ with respect to $(\btheta
,\bl)$ follows from the measurability of the coefficients $H^n$, $\alpha
_T^n$ and $f^n$ w.r.t. $(\btheta,\bl)$ (see Appendix C in \cite
{khalim11}). The property $(\Pc_n)$ is then proved.

$\bullet$ Fix\vspace*{1pt} $k \in\{0,\ldots,n-1\}$, and suppose that
$(\Pc_{k+1})$ is true, and denote by $(Y^{k+1},Z^{k+1}) \in\Sb
_c^\infty(\Delta_{k+1},E^{k+1}) \times\Lb^2_W(\Delta
_{k+1},E^{k+1})$ a solution to ${(E_{k+1})}$. Since the indexed $\F
$-adapted process
$Y^{k+1}$ is continuous, it is actually $\F$-predictable, and so
$Y^{k+1} \in\Pc_\F(\Delta_{k+1},E^{k+1})$. This implies that the
map $(t,\omega,\btheta_k,\bl_{k+1}) \rightarrow Y_t^{k+1}(\omega
,\btheta_k,t,\bl_{k+1})$ is
$\Pc(\F)\otimes\Bc(\Delta_k)\otimes\Bc(E^{k+1})$-measurable. The
generator $f^{k}$ is thus well defined in \eqref{deffk} as an indexed process
in $\Pc_\F(\R,\R^m,\Delta_k,E^k)$, and we shall prove that (\ref{Pck})
holds true by proceeding in four steps, in order to overcome the
technical difficulties coming from the exponential term
in $U(y)$ together with the quadratic condition in $z$ for
$f^k$.\vspace*{10pt}


\textit{Step} 1: \textit{Approximating sequence}. We truncate
the term $U(y) = -e^{-py}$ when $y$ goes to $-\infty$, as well as the
infimum, by considering the truncated generator
\begin{eqnarray*}
&&f_N^{k}(t,y,z,\btheta_k,\bl_k)\\[-0.5pt]
&&\qquad= - \lambda_t^{k}(\btheta_k,\bl_k)  z - \frac{1}{2p} |\lambda
_t^{k}(\btheta_k,\bl_k)|^2 \\[-0.5pt]
&&\qquad\quad{}+ \inf_{\pi\in A^{k},|(\sigma_t^k)'\pi|\leq N} \biggl\{ \frac{p}{2}
\biggl| z + \frac{1}{p} \lambda_t^{k}(\btheta_k,\bl_k) - \sigma
_t^{k}(\btheta_k,\bl_k)'\pi\biggr|^2 \\[-0.5pt]
&&\qquad\quad\hspace*{83.5pt}{}+ \frac{1}{p} U(\max(-N,y))\\[-0.5pt]
&&\qquad\quad\hspace*{96.7pt}\times\int_E U\bigl(\pi\gamma_t^k(\btheta
_k,\bl
_k,\ell)\\[-0.5pt]
&&\qquad\quad\hspace*{132pt}{}
- Y_t^{k+1}(\btheta_k,t,\bl_k,\ell)\bigr) \eta_{{k+1}}(\bl_k,d\ell)
\biggr\}
\end{eqnarray*}
and introduce the corresponding family of approximated BSDEs with
terminal data
$\tilde H_T^{k}$ and generator $f_N^{k}$,
%
\setcounter{equation}{4}
\renewcommand{\theequation}{\arabic{section}.\arabic{equation}}
\begin{eqnarray} \label{BSDEYkN}
Y_t^{k,N}(\btheta_k,\bl_k) &=& H_T^{k}(\btheta_k,\bl_k) +\frac{1}{p}\ln
\alpha_T^k(\btheta_k,\bl_k) \nonumber\\
&&{} + \int_t^T f_N^{k}(r,Y_r^k,Z_r^{k,N},\btheta_k,\bl_k) \,dr\\
&&{}- \int_t^T Z_r^{k,N}\,dW_r,\qquad \theta_{k} \leq t \leq T. \nonumber
\end{eqnarray}
Under~\ref{HB}(i), there exists a constant $C$ such that for all
\mbox{$(\btheta_k,\bl_k) \in\Delta_k(T)\times E^k$},
%
\begin{eqnarray}\label{BIfIn}
f_N^{k}(t,y,z,\btheta_k,\bl_k) & \geq& - \lambda_t^{k}(\btheta_k,\bl
_k)  z - \frac{1}{2p}
|\lambda_t^{k}(\btheta_k,\bl_k)|^2\nonumber\\[-8pt]\\[-8pt]
& \geq& - C(1+|z|)\nonumber
\end{eqnarray}
for all $(t,y,z)\in[\theta_k,T]\times\R\times\R^m$. Moreover, since $0
\in A^k$, and the process $Y^{k+1}$ is essentially bounded,
there exists some positive constant $C_N$ (depending on $N$) s.t. for
all $(\btheta_k,\bl_k) \in\Delta_k(T)\times E^k$,
%
\begin{eqnarray}\label{BSfIn}
f_N^{k}(t,y,z,\btheta_k,\bl_k) & \leq& - \lambda_t^{k}(\btheta_k,\bl
_k)  z - \frac{1}{2p} | \lambda_t^{k}(\btheta_k,\bl_k)|^2
\nonumber\\
&&{}+ {p\over2 }\biggl| z + \frac{1}{p} \lambda_t^{k}(\btheta_k,\bl_k)
\biggr|^2 + C_N \\
& \leq& C_N (1 + |z|^2),\nonumber
\end{eqnarray}
under~\ref{HB}(i). Hence, for any given $(\btheta_k,\bl_k) \in
\Delta_k(T)\times E^k$, we can apply
Theorem~2.3 in~\cite{kob00}, and obtain the existence of a solution
$(Y^{k,N}(\btheta_k,\bl_k),\break Z^{k,N}(\btheta_k$, $\bl_k)) \in\Sb_c^\infty
[\theta_k,T]\times\Lb_W^2[\theta_{k},T]$ to \eqref
{BSDEYkN}. The measurability of $(Y^{k,N},\break Z^{k,N})$ w.r.t. its arguments
$(\btheta_k,\bl_k)$ follows from the measurability of $H_T^k$, $\alpha
_T^k$, $f^{k}_N$ w.r.t. $(\btheta_k,\bl_k)$. In the next steps, we
prove the convergence of the sequence $(Y^{k,N},Z^{k,N})_{N}$ to a
solution of (\ref{EK}).\vspace*{10pt}

\textit{Step} 2: \textit{Lower bound for the approximating
sequence}. Define the generator function $\underline{f}^k$ by
\[
\underline{f}^k(t,z,\btheta_{k},\bl_{k}) = - \lambda_t^{k}(\btheta_k,\bl
_k)  z - \frac{1}{2p} |\lambda_t^{k}(\btheta_k,\bl_k)|^2.
\]
Under\vspace*{1pt}~\ref{HB}(i), and for fixed $(\btheta_k,\bl_k) \in\Delta
_k(T)\times E^k$, the function $\underline f^k(\cdot,\btheta_{k},\bl_{k})$
satisfies the usual Lipschitz (and a fortiori quadratic growth)
condition in $z$, which implies from Theorem 2.3 in~\cite{kob00} that
there exists
$( \underline Y^k(\btheta_{k},\bl_{k}),\break \underline Z^k(\btheta_{k},\bl
_{k})) \in\Sb_c^\infty[\theta_k,T]\times\Lb_W^2[\theta_k,T]$ solution
to the
BSDE with terminal data
$H_T^{k}(\btheta_k,\bl_k) + \frac{1}{p}\ln\alpha_T^k(\btheta_k,\bl
_k)$, and generator $\underline{f}^k(\cdot, \cdot,\btheta_{k},\bl_{k})$.
The solution $( \underline Y^k, \underline Z^k)$ is measurable w.r.t.
the arguments $(\btheta_k,\bl_k)$, and from the uniform boundedness
condition in~\ref{HB}, and Proposition 2.1 in~\cite{kob00}, we
deduce that $( \underline Y^k, \underline Z^k) \in\Sb_c^\infty(\Delta
_k,E^k)\times\Lb_W^2(\Delta_k,E^k)$. Moreover, we
easily see under~\ref{HB}(i)
that for any $(\btheta_k,\bl_k) \in\Delta_k(T)\times E^k$,
$\underline f^k(\cdot,\btheta_{k},\bl_{k})$ satisfy Assumptions (H2) and
(H3) of~\cite{kob00}.
Since $\underline f^k(\cdot,\btheta_{k},\bl_{k}) \leq f^{k}_N(\cdot
,\btheta
_{k},\bl_{k})$, we can apply comparison Theorem 2.6 in~\cite{kob00} to
get the inequality
%
\begin{equation} \label{minYn}
Y^{k,N}_{t}(\btheta_{k},\bl_{k}) \geq\underline Y^k_{t}(\btheta_{k},\bl
_{k}),\qquad \theta_k\leq t \leq T \qquad\mbox{a.s.}
\end{equation}
for all $N$, and $(\btheta_k,\bl_k) \in\Delta_k(T)\times E^k$.
Since $\underline Y^k \in\Sb_c^\infty(\Delta_k,E^k)$, this implies that
$Y^{k,N}$ is uniformly lower bounded, and thus by \eqref{BSDEYkN}, we
see\vspace*{1pt} that for $N$ large enough,
$(Y^{k,N},Z^{k,N})$ satisfies the indexed BSDE with terminal data
$\tilde H_T^{k}$, and with a generator $\tilde f_N^k$ where one can
remove in
$f_N^k$ the truncation in $-N$ for $U(y)$, that is,
%
\begin{eqnarray} \label{BSDEYkNtilde}
Y_t^{k,N}(\btheta_k,\bl_k) &=& H_T^{k}(\btheta_k,\bl_k) +\frac{1}{p}\ln
\alpha_T^k(\btheta_k,\bl_k)\nonumber\\[-2pt]
&&{} + \int_t^T \tilde f_N^{k}(r,Y_r^k,Z_r^{k,N},\btheta_k,\bl
_k) \,dr\\[-2pt]
&&{} - \int_t^T Z_r^{k,N}\,dW_r,\qquad \theta_{k} \leq t \leq T,
\nonumber
\end{eqnarray}
with
\begin{eqnarray*}
&&\tilde f_N^{k}(t,y,z,\btheta_k,\bl_k)\\[-2pt]
&&\qquad= - \lambda_t^{k}(\btheta_k,\bl_k)  z - \frac{1}{2p} |\lambda
_t^{k}(\btheta_k,\bl_k)|^2 \\[-2pt]
&&\qquad\quad{} + \inf_{\pi\in A^{k},|(\sigma_t^k)'\pi|\leq N} \biggl\{ \frac
{p}{2} \biggl| z + \frac{1}{p} \lambda_t^{k}(\btheta_k,\bl_k) - \sigma
_t^{k}(\btheta_k,\bl_k)'\pi\biggr|^2 \\[-2pt]
&&\qquad\quad\hspace*{84.3pt}{} + \frac{1}{p} U(y) \int_E U\bigl(\pi\gamma_t^k(\btheta_k,\bl
_k,\ell)\\[-2pt]
&&\qquad\quad\hspace*{156.6pt}{} - Y_t^{k+1}(\btheta_k,t,\bl_k,\ell)\bigr) \eta_{{k+1}}(\bl_k,d\ell)
\biggr\}.
\end{eqnarray*}

\textit{Step} 3: \textit{Monotonicity and uniform estimate of
the approximating sequence.} We cannot apply directly a comparison theorem
for $Y^{k,N}$ for the quadratic generators $\tilde f_N^k$, since the
derivative of $\tilde f_N^k$, with respect to $z$, is not of linear
growth in $z$, as
requested in Assumption (H2) in~\cite{kob00}. We then make an
exponential change of variable by defining for any
$(\btheta_k,\bl_k) \in\Delta_k(T)\times E^k$, the pair of
processes $(\dot Y ^{k,N}(\btheta_k,\bl_k), \dot Z ^{k,N}(\btheta_k,\bl
_k)) \in\Sb_c^\infty[\theta_k,T]\times\Lb_W^2[\theta_k,T]$ by
\[
\dot Y^{k,N}_{t} (\btheta_{k},\bl_{k}) = \exp(p Y^{k,N}_{t}
(\btheta_{k},\bl_{k}))
\]
and
\[
\dot Z^{k,N}_{t} (\btheta_{k},\bl_{k}) = p\dot Y^{k,N}_{t}
(\btheta_{k},\bl_{k}) Z^{k,N}_{t} (\btheta_{k},\bl_{k}).
\]
A straightforward It\^o formula on \eqref{BSDEYkNtilde} shows that
$(\dot Y^{k,N}(\btheta_{k},\bl_{k}), \dot Z^{k,N}(\btheta_{k},\bl_{k}))$
is solution to the BSDE
\begin{eqnarray*}
\dot Y^{k,N}_{t}(\btheta_{k},\bl_{k}) & = & \alpha_{T}^k(\btheta_{k},\bl
_{k}) \exp(pH_{T}^k(\btheta_{k},\bl_{k}))\\[-2pt]
&&{} + \int_{t}^T\dot f_N^{k}(r,\dot Y^{k,N}_{r},\dot
Z^{k,N}_{r},\btheta_{k},\bl_{k})\,dr\\[-2pt]
&&{} -
\int_{t}^T\dot Z^{k,N}_{r}\,dW_{r}, \qquad\theta_{k} \leq t \leq T,
\end{eqnarray*}
where the generator $\dot f_N^{k}$ is defined by
\begin{eqnarray*}
&&
\dot f_N^{k}(t,y,z,\btheta_{k},\bl_{k})\\[-2pt]
&&\qquad =
\inf_{\pi\in A^k,|(\sigma_t^k)'\pi|\leq N} \biggl\{ \frac{1}{2}p^2y|\sigma
^k_{t}(\btheta_{k},\bl_{k})'\pi|^2\\[-2pt]
&&\qquad\hspace*{82.1pt}{}-p\bigl(\lambda^k_{t}(\btheta_{k},\bl_{k})y+z\bigr) \sigma^k_{t}(\btheta
_{k},\bl
_{k})'\pi\\[-2pt]
&&\hspace*{70.5pt}\qquad\quad{} - \int_{E} U\bigl(\pi \gamma_{t}^{k}(\btheta_{k},\bl
_{k},\ell)\\[-2pt]
&&\hspace*{142pt}{} - Y^{k+1}_{t}(\btheta_{k},t,\bl_{k},\ell)\bigr)\eta_{{k+1}}(\bl
_k,d \ell)\biggr\}.
\end{eqnarray*}
Fix $(\btheta_k,\bl_k) \in\Delta_k\times E^k$.
Denote by $\dot g_N^{k}(\pi,t,y,z,\btheta_{k},\bl_{k})$ the function
inside the infimum defining $\dot f_N^{k}$, that is, $\dot f_N^k(\cdot
) = \inf
_{\pi\in A^k, |(\sigma_t^k)'\pi|\leq N} \dot g_N^k(\pi,\cdot)$. Then,
for all $(t,y,y',z,z',\btheta_k,\bl_k) \in[\theta_k,T]\times\R^2\times
(\R^m)^2\times\Delta_k\times E^k$, we have
\begin{eqnarray*}
& & | \dot f_N^{k}(t,y,z,\btheta_{k},\bl_{k}) - \dot
f_N^{k}(t,y',z',\btheta_{k},\bl_{k}) | \\
&&\qquad \leq \sup_{\pi\in A^k,|(\sigma_t^k)'\pi|\leq N} | \dot
g_N^{k}(\pi,t,y,z,\btheta_{k},\bl_{k}) - \dot g_N^{k}(\pi
,t,y',z',\btheta_{k},\bl_{k}) | \\
&&\qquad \leq \biggl(\frac{1}{2}p^2 N + pN |\lambda^k_{t}(\btheta_{k},\bl
_{k})| \biggr)|y-y'| + pN |z-z'|.
\end{eqnarray*}
Under~\ref{HB}(i), we then see that $\dot f_N^{k}$ satisfies the
standard Lipschitz condition in $(y,z)$, uniformly in $(t,\omega)$.
Since the sequence $(\dot f_N^k)_N$ is noninceasing, that is, $\dot
f_{N+1}^k \leq\dot f_N^k$, we obtain\vspace*{1pt} by standard
comparison principle for BSDE that $\dot Y ^{k,N+1} \leq\dot Y ^{k,N}$,
and so
%
\begin{equation}\label{decrYn}\qquad
Y^{k,N+1}_{t}(\btheta_k,\bl_k) \leq Y^{k,N}
_{t}(\btheta_k,\bl_k),\qquad
\theta_k\leq t \leq T \qquad\mbox{a.s. } \forall N \in\mathbb{N}
\end{equation}
for all $(\btheta_k,\bl_k) \in\Delta_k\times E^k$. From the
quadratic condition in $z$ for $f_0^k$ in \eqref{BIfIn} and \eqref{BSfIn},
uniformly in
$(\btheta_k,\bl_k)$, and the a priori estimate of Proposition 2.1 in
\cite{kob00}, we deduce under~\ref{HB}(ii) that $Y^{k,0} \in\Sb
_c^\infty(\Delta_k,E^k)$. Together with \eqref{minYn} and \eqref
{decrYn}, this implies that there exists a positive constant $M$ such that
%
\begin{equation} \label{UBYn}
\|Y^{k,N}\|_{\Sb_c^\infty(\Delta_k,E^k)} \leq M\qquad\forall N \in
\mathbb{N}.
\end{equation}

\vspace*{5pt}

\textit{Step} 4: \textit{Convergence of the approximating sequence.} By
using \eqref{UBYn} in \eqref{BSDEYkN} [or~\eqref{BSDEYkNtilde}], we
see that $(Y^{k,N},Z^{k,N})$ satisfies the indexed BSDE with terminal
data~$\tilde H_T^k$, and with generator $\hat f^{k}_N$ given by
\begin{eqnarray*}
&&
\hat f_N^{k}(t,y,z,\btheta_k,\bl_k) \\
&&\qquad= - \lambda_t^{k}(\btheta_k,\bl
_k)  z - \frac{1}{2p} |\lambda_t^{k}(\btheta_k,\bl_k)|^2 \\
&&\qquad\quad{} + \inf_{\pi\in A^{k},|(\sigma_t^k)'\pi|\leq N} \biggl\{ \frac
{p}{2} \biggl| z + \frac{1}{p} \lambda_t^{k}(\btheta_k,\bl_k) - \sigma
_t^{k}(\btheta_k,\bl_k)'\pi\biggr|^2 \\
&&\hspace*{84pt}\qquad\quad{} + \frac{1}{p} U\bigl((-M)\vee y\bigr)\\
&&\qquad\quad\hspace*{95pt}{}\times \int_E U\bigl(\pi\gamma
_t^k(\btheta_k,\bl_k,\ell)\\
&&\qquad\quad\hspace*{133pt}{} - Y_t^{k+1}(\btheta_k,t,\bl_k,\ell)\bigr) \eta_{{k+1}}(\bl_k,d\ell)
\biggr\}.
\end{eqnarray*}
By the same arguments as for the generator $f_N^{k}$, there exists a
constant $C_M$ such that
\[
|\hat f_N^{k}(t,y,z,\btheta_k,\bl_k)| \leq C_M(1+|z|^2)
\]
for all $N \in\mathbb{N}$, $(t,y,z)\in[0,T]\times\R\times\R^m$,
$(\btheta_k,\bl_k) \in\Delta_k\times E^k$. Let us check that
the nonincreasing sequence $(\hat f_N^{k})_N$ converges uniformly on
compact sets of $(t,y,z) \in[0,T]\times\R\times\R^m$ to $\hat
f^{k}$ defined by
\begin{eqnarray*}
&&
\hat f^{k}(t,y,z,\btheta_k,\bl_k) \\
&&\qquad= - \lambda_t^{k}(\btheta_k,\bl
_k)  z - \frac{1}{2p} |\lambda_t^{k}(\btheta_k,\bl_k)|^2 \\
&&\qquad\quad{} + \inf_{\pi\in A^{k}} \biggl\{ \frac{p}{2} \biggl| z + \frac
{1}{p} \lambda_t^{k}(\btheta_k,\bl_k) - \sigma_t^{k}(\btheta_k,\bl
_k)'\pi\biggr|^2 \\
&&\qquad\quad\hspace*{39pt}{} + \frac{1}{p} U\bigl((-M)\vee y\bigr) \int_E U\bigl(\pi \gamma
_t^k(\btheta_k,\bl_k,\ell)\\
&&\qquad\quad\hspace*{151pt}{}- Y_t^{k+1}(\btheta_k,t,\bl_k,\ell)\bigr) \eta_{{k+1}}(\bl_k,d\ell)
\biggr\}.
\end{eqnarray*}
Indeed, notice that in the definition of $\hat f^k$, one may restrict
in the infimum over $\pi$ in $A^k$ s.t. the function $\hat g^k(\pi
,\cdot)$
inside the infimum bracket, that is,
\begin{eqnarray*}
&&
\hat g^k(\pi,t,y,z,\btheta_k,\bl_k) \\
&&\qquad= \frac{p}{2} \biggl| z + \frac
{1}{p} \lambda_t^{k}(\btheta_k,\bl_k) - \sigma_t^{k}(\btheta_k,\bl
_k)'\pi\biggr|^2 \\
&&\qquad\quad{} + \frac{1}{p} U\bigl((-M)\vee y\bigr) \int_E U\bigl(\pi \gamma_t^k(\btheta
_k,\bl_k,\ell)\\
&&\qquad\quad\hspace*{111.5pt}{} - Y_t^{k+1}(\btheta_k,t,\bl_k,\ell)\bigr) \eta_{{k+1}}(\bl_k,d\ell)
\end{eqnarray*}
is smaller than $\hat g^k(\pi,\cdot)$ for $\pi= 0$. In other words, we have
\begin{eqnarray*}
\hat f^{k}(t,y,z,\btheta_k,\bl_k) &=& - \lambda_t^{k}(\btheta_k,\bl_k)
z
- \frac{1}{2p} |\lambda_t^{k}(\btheta_k,\bl_k)|^2 \\
&&{}+
\inf_{\pi\in A^k \cap\mathcal{K}(t,y,z,\btheta_k,\bl_k)} \hat g^k(\pi
,t,y,z,\btheta_k,\bl_k),
\end{eqnarray*}
where
\[
\mathcal{K}(t,y,z,\btheta_k,\bl_k) = \{ \pi\in\R^d\dvtx \hat g^k(\pi
,t,y,z,\btheta_k,\bl_k) \leq\hat g^k(0,t,y,z,\btheta_k,\bl_k)
\}.
\]
Since $U$ is nonpositive, $Y^{k+1}$ is essentially bounded, and under
\ref{HB}(i), there exists some positive constant $C$ such that
%
\begin{eqnarray} \label{pileqz}
&&
\mathcal{K}(t,y,z,\btheta_k,\bl_k) \nonumber\\
&&\qquad \subset \biggl\{ \pi\in\R^d\dvtx
\biggl| z + \frac{1}{p} \lambda_t^{k}(\btheta_k,\bl_k) - \sigma
_t^{k}(\btheta_k,\bl_k)'\pi\biggr| \leq
\biggl| z + \frac{1}{p} \lambda_t^{k}(\btheta_k,\bl_k) \biggr| + C \biggr\}
\hspace*{-10pt}\\
&&\qquad \subset \{ \pi\in\R^{d}\dvtx |\sigma_t^{k}(\btheta_k,\bl
_k)'\pi| \leq C(|z|+1)\}\nonumber
\end{eqnarray}
for all $(t,y,z,\btheta_k,\bl_k) \in[0,T]\times\R\times\R^m\times
\Delta_k\times E^k$.
This shows that on any compact of $(t,y,z) \in[0,T]\times\R\times\R
^m$, we have $\mathcal{K}(t,y,z,\btheta_k,\bl_k) \subset\{ \pi\dvtx
|(\sigma_t^k)'\pi| \leq N \}$ for $N$ large enough,
and so
$\hat f_N^k = \hat f^k$, which obviously implies
the convergence of $(\hat f_N^{k})_N$ to $\hat f^k$ locally uniformly
on $(t,y,z) \in[0,T]\times\R\times\R^m$.
We can then apply Proposition 2.4 in~\cite{kob00}, which states that
the sequence $(Y^{k,N}(\btheta_{k},\bl_{k}),Z^{k,N}(\btheta_{k},\bl
_{k}))_{N}$ converges in $\Sb_c^\infty[\theta_{k},T]\times\Lb_W^2[\theta
_{k},T]$ to $(Y^{k}(\btheta_{k}$, $\bl_{k}),Z^{k}(\btheta_{k},\bl_{k}))$
solution to the BSDE with terminal data $\tilde H_T^k$, and generator
$\hat f^{k}$. The indexed processes $(Y^k,Z^k)$ inherit from
$(Y^{k,N},Z^{k,N})$ the measurability in
the arguments $(\btheta_k,\bl_k) \in\Delta_k\times E^k$. Moreover,
from \eqref{UBYn}, we see that $Y^k$ also satisfies the estimate
\[
\|Y^{k}\|_{\Sb_c^\infty(\Delta_k,E^k)} \leq M.
\]
Hence, this implies that one can remove the truncation term $-M$ in the
BSDE with generator $\hat f^k$ satisfied by $(Y^k,Z^k)$. Therefore,
$(Y^k,Z^k) \in\Sb_c^\infty(\Delta_k,\break E^k)\times\Lb_W^2(\Delta
_k, E^k)$ is solution to (\ref{EK}), which ends the induction proof.
\end{pf}

\subsection{BSDE characterization by verification theorem}\label{sec4.2}

In this section, we show that a solution $(Y^k)_k$ to the recursive
system indexed BSDEs actually provides the solution to the optimal
investment problem in terms of the value functions $V^k$, $k = 0,\ldots
,n$, in \eqref{defVkdyn}. As a byproduct, we get the uniqueness
of this system of BSDEs and a description of an optimal strategy by
means of the solution to these BSDEs.
%
\begin{Thm}
The value functions $V^k$, $k = 0,\ldots,n$, defined in \eqref
{defVndyn}, \eqref{defVkdyn}, from the decomposition of the optimal
investment problem \eqref{defV0}, are given by
%
\begin{equation} \label{VkYk}
V_t^k(x,\btheta_k,\bl_k,\nu^k) = U\bigl(X_t^{k,x} - Y_{t}^{k}(\btheta
_k,\bl_k)\bigr),\qquad \theta_k \leq t \leq T,
\end{equation}
for all $x \in\R$, $(\btheta_k,\bl_k) \in\Delta_k\times E^k$,
$\nu^k \in\Ac_\F^k$, where
\[
(Y^k,Z^k)_{k=0,\ldots,n} \in\prod_{k=0}^n
\Sb_c^\infty(\Delta_k,E^k) \times\Lb^2_W(\Delta
_k,E^k)
\]
is the solution to the recursive system of indexed BSDEs (\ref{EK}),
$k = 0,\ldots,n$. Here, $X^{k,x}$ denotes the wealth process in\vadjust{\goodbreak}
\eqref{dynXI} controlled by $\nu^k$, and starting from $x$ and~$\theta
_k$. Moreover, there exists an optimal trading strategy $\hat\pi= (\hat
\pi^k)_{k=0,\ldots,n} \in\Ac_\G= (\Ac_\F^k)_{k=0,\ldots
,n}$ described by
%
\begin{eqnarray} \label{formpik}
&&\hat\pi_t^k(\btheta_k,\bl_k) \nonumber\\[-2pt]
&&\qquad \in \argmin_{\pi\in A^k} \biggl\{ \frac
{p}{2} \biggl| Z_t^k(\btheta_k,\bl_k) + \frac{1}{p} \lambda_t^{k}(\btheta
_k,\bl_k) - \sigma_t^{k}(\btheta_k,\bl_k)'\pi\biggr|^2
\nonumber\\[-9pt]\\[-9pt]
&&\qquad\quad\hspace*{37.4pt}{} + \frac{1}{p} U(Y_t^k(\btheta_k,\bl_k)) \int_E U
\bigl(\pi\gamma_t^k(\btheta_k,\bl_k,\ell)\nonumber\\[-2pt]
&&\qquad\quad\hspace*{150pt}{} - Y_t^{k+1}(\btheta_k,t,\bl_k,\ell)\bigr) \eta_{{k+1}}(\bl_k,d\ell)
\biggr\} \nonumber
\end{eqnarray}
for $k = 0,\ldots,n-1$, $(\btheta_k,\bl_k) \in\Delta
_k(T)\times E^k$, $t \in[\theta_k,T]$, a.s., and
\[
\hat\pi_t^n(\btheta,\bl) \in\argmin_{\pi\in A^n}
\biggl| Z_t^n(\btheta,\bl) + \frac{1}{p} \lambda_t^{n}(\btheta,\bl) -
\sigma_t^{n}(\btheta,\bl)'\pi\biggr|^2
\]
for $k = n$, $(\btheta,\bl) \in\Delta_n(T)\times E^n$, $t \in
[\theta_n,T]$, a.s.
\end{Thm}
\begin{pf}
\textit{Step} 1:
\textit{We first prove that for all $k = 0,\ldots,n$, $\nu^k \in\Ac_\F^k$,
$U(X^{k,x}-Y^k(\btheta_k,\bl_k)) \geq V^k(x,\btheta_k,\bl_k,\nu^k)$.}
Let $(Y^k,Z^k)_{k=0,\ldots,n} \in\prod_{k=0}^n \Sb_c^\infty(\Delta
_k,E^k) \times\Lb^2_W(\Delta
_k,E^k)$ be a solution to the system of BSDEs (\ref{EK}), $k = 0,\ldots,n$.
For any $x \in\R$, $(\btheta_k,\bl_k) \in\Delta_k(T)\times
E^k$, $\nu^k \in\Ac_\F^k$, we apply It\^o's formula to the process
\begin{eqnarray*}
\xi_t^k(x,\btheta_k,\bl_k,\nu^k) &:=& U\bigl(X_t^{k,x} -
Y_{t}^{k}(\btheta_k,\bl_k)\bigr) \\[-2pt]
& &{} + \int_{\theta_k}^t \int_E U\bigl(X_s^{k,x} + \nu_s^k \gamma
_s^k(\btheta_k,\bl_k,\ell)\\[-2pt]
&&\hspace*{65pt}{} - Y_s^{k+1}(\btheta_k,s,\bl_k,\ell) \bigr)
\eta_{{k+1}}(\bl_k,d\ell)\,ds
\end{eqnarray*}
for $k = 0,\ldots,n$, and $\xi_t^n(x,\btheta_k,\bl_n,\nu^n):=
U(X_t^{n,x} - Y_{t}^{n}(\btheta_n,\bl_n))$,
for $k = n$, and $\theta_k\leq t\leq T$.
From the dynamics of $X^{k,x}$ and $Y^k$, we immediately get
\begin{eqnarray*}
& & d \xi_t^k(x,\btheta_k,\bl_k,\nu^k) \\[-2pt]
&&\qquad= - U\bigl(X_t^{k,x} - Y_{t}^{k}(\btheta_k,\bl_k)\bigr)
\bigl[ \bigl(f_t^k(t,Y_t^k,Z_t^k,\btheta_k,\bl_k)\\[-2pt]
&&\qquad\quad\hspace*{115pt}{} - g_t^k(\nu
_t^k,t,Y_t^k,Z_t^k,\btheta_k,\bl_k) \bigr) \,dt \\[-2pt]
&&\hspace*{123.5pt}\qquad\quad{} + \bigl( \sigma_t^k(\btheta_k,\bl_k)'\nu_t^k -
Z_t^k\bigr)\,dW_t \bigr],
\end{eqnarray*}
where
\begin{eqnarray*}
&&
g_t^k(\pi,t,y,z,\btheta_k,\bl_k)\\[-2pt]
&&\qquad= \frac{p}{2} | z - \sigma
_t^{k}(\btheta_k,\bl_k)'\pi|^2 - b_t^{k}(\btheta_k,\bl_k)'\pi\\[-2pt]
&&\qquad\quad{} + \frac{1}{p} U(y) \int_E U\bigl(\pi\gamma_t^k(\btheta
_k,\bl_k,\ell)\\[-2pt]
&&\qquad\quad\hspace*{71.6pt}{} - Y_t^{k+1}(\btheta_k,t,\bl_k,\ell)\bigr)
\eta_{{k+1}}(\bl_k,d\ell)
\end{eqnarray*}
for $k = 0,\ldots,n-1$, and $g_t^n(\pi,t,y,z,\btheta_n,\bl_n) = \frac
{p}{2} | z - \sigma_t^{n}(\btheta_n,\bl_n)'\pi|^2 -
b_t^{n}(\btheta_n$, $\bl_n)'\pi$ for $k = n$.
Since, by construction, $f_t^k(t,y,z,\btheta_k,\bl_k) = \inf_{\pi
\in A_k} g_t^k(\pi,t,y,z$, $\btheta_k,\bl_k)$, and recalling that $U$ is
nonpositive, this implies that the process $\{\xi_t^k(x,\btheta_k$, $\bl
_k,\nu^k)$, $\theta_k\leq t\leq T\}$, is a local supermartingale. By
considering a localizing
$\F$-stopping times sequence $(\rho_n)_n$ valued in $[\theta_k,T]$ for
$\xi^k$, we have the inequality
\[
\Ee[ \xi_{s\wedge\rho_n}^k (x,\btheta_k,\bl_k,\nu^k) | \Fc_t
] \leq
\xi_{t\wedge\rho_n}^k(x,\btheta_k,\bl_k,\nu^k), \qquad\theta_k \leq t
\leq s \leq T.
\]
Now, by Definition~\ref{defadmi} of the admissibility condition for
$\nu ^k$, and since the processes $Y^k$, $Y^{k+1}$ are essentially
bounded, the sequence
$(\xi_{s\wedge\rho_n}^k(x,\btheta_k,\bl_k,\break\nu^k))_n$ is uniformly
integrable for any $s \in[\theta_k,T]$, and by the dominated
convergence theorem, we obtain the supermartingale property of
$\xi^k(x,\btheta_k,\bl_k,\nu^k)$. Therefore, by writing the
supermartingale\vspace*{1pt} property between $t$ and $T$, and
recalling that $Y_T^k = H_T^k + \frac{1}{p}\ln\alpha_T^k$, we obtain
the inequalities
%
\begin{equation} \label{Unsuper}
U\bigl(X_t^{n,x} - Y_{t}^{n}(\btheta,\bl)\bigr)
\geq \Ee\bigl[
U\bigl(X_T^{n,x} - H_T^n(\btheta,\bl)\bigr) \alpha_T(\btheta,\bl) | \Fc_t
\bigr],
\end{equation}
\begin{eqnarray}
\label{Uksuper}
&&U\bigl(X_t^{k,x} - Y_{t}^{k}(\btheta_k,\bl_k)\bigr)
\nonumber\\
&&\qquad \geq \Ee\biggl[
U\bigl(X_T^{k,x} - H_T^k(\btheta_k,\bl_k)\bigr) \alpha_T^k(\btheta_k,\bl_k)
\nonumber\\[-8pt]\\[-8pt]
&&\qquad\quad\hspace*{10.3pt}{} + \int_{t}^T \int_E U\bigl(X_s^{k,x} + \nu_s^k
\gamma_s^k(\ell)\nonumber\\
&&\qquad\quad\hspace*{66.2pt}{} -
Y_s^{k+1}(\btheta_k,s,\bl_k,\ell) \bigr) \eta_{{k+1}}(\bl_k,d\ell)\,ds
\Big| \Fc_t \biggr], \nonumber
\end{eqnarray}
which hold true for any $\nu^k \in\Ac_\F^k$, $k =
0,\ldots,n$.\vspace*{10pt}

\textit{Step} 2: \textit{The process $\int_{\theta_k}^{ \cdot}
Z^{k}_{s}(\btheta_{k},\bl_{k})\,dW_{s}$ is a BMO-martingale, for any $k
= 0,\ldots,n$,
$(\btheta_k,\bl_k) \in\Delta_k(T)\times E^k$}. By applying It\^o's
formula to the process
$\exp(qY_t^{k}(\btheta_{k},\bl_{k}))$ with $q>p$ between any stopping
time $\tau$ valued in $[\theta_{k}, T]$ and~$T$, and recalling
the terminal data $Y_T^k = \tilde H_T^k = H_T^k+\frac{1}{p}\ln
\alpha_T^k$, we get
%
\begin{eqnarray} \label{decompBMO1}\quad
& & \frac12 q(q-p)\mathbb{E}\biggl[\int_{\tau}^T\exp
(qY_{t}^{k}(\btheta_{k},\bl_{k}))|Z_{t}^{k}(\btheta_{k},\bl
_{k})|^2\,dt\Big|\Fc_{\tau}\biggr] \nonumber\\
&&\qquad= q\mathbb{E}\biggl[\int_{\tau}^T\exp(qY_{t}^{k}(\btheta_{k},\bl
_{k}))\biggl(f^k(t,Y_{t}^{k},Z_{t}^{k},\btheta_k,\bl_{k})
-\frac{p}{2} |Z_{t}^{k}|^2 \biggr) \,dt \Big|\Fc_{\tau} \biggr] \\
&&\qquad\quad{} + \mathbb{E}[\exp(q \tilde H_{T}^{k}(\btheta_{k},\bl
_{k}))
- \exp(qY_{\tau}^{k}(\btheta_{k},\bl_{k})) |\Fc_{\tau}
].\nonumber
\end{eqnarray}
By definition of $f^k$ in \eqref{deffk}, and since $Y^{k+1} \in\Sb
_c^\infty(\Delta_{k+1},E^{k+1})$,
there exists a constant $C$ such that for all $(t,y,z)\in[0,T]\times\R
\times\R^d$,
\[
f^k(t,y,z,\btheta_{k}, \bl_{k}) \leq\frac{p}{2} |z|^2 - C U(y).
\]
Combining this last inequality with \eqref{decompBMO1}, we get
\begin{eqnarray*}
& & \frac12 q(q-p)\mathbb{E}\biggl[\int_{\tau}^T\exp
(qY_{t}^{k}(\btheta_{k},\bl_{k}))|Z_{t}^{k}(\btheta_{k},\bl_{k})|^2
\,dt \Big|\Fc_{\tau}\biggr] \\
&&\qquad \leq qC\mathbb{E}\biggl[\int_{\tau}^T \exp\bigl((q-p)Y_{t}^{k}(\btheta
_{k},\bl_{k})\bigr) \,dt\Big|\Fc_{\tau} \biggr]\\
&&\qquad\quad{}+\mathbb{E}\bigl[ e^{q\tilde H_{T}^{k}(\btheta_{k},\bl_{k})} -
e^{qY_{\tau}^{k}(\btheta_{k},\bl_{k})} |\Fc_{\tau} \bigr].
\end{eqnarray*}
Under~\ref{HB}(ii), and since $Y^{k} \in\Sb_c^\infty(\Delta
_{k},E^{k})$, this shows that there exists a constant $C$ s.t.
\[
\mathbb{E}\biggl[\int_{\tau}^T |Z_{t}^{k}(\btheta_{k},\bl_{k})|^2 \,dt
\Big|\Fc_{\tau}\biggr] \leq C\qquad \mbox{for any stopping time } \tau
\mbox{ valued in }[\theta_k,T],
\]
which is the required BMO-property.\vspace*{10pt}

\textit{Step} 3: \textit{Admissibility of $\hat\pi^k$}. Let us
consider the functions $\hat g^k$, $k = 0,\ldots,n$, defined by
\begin{eqnarray*}
&&
\hat g^k(\pi,t,\omega,\btheta_k,\bl_k) \\
&&\qquad=
\frac{p}{2} \biggl| Z^k_t(\btheta_k,\bl_k) + \frac{1}{p} \lambda
_t^{k}(\btheta_k,\bl_k) - \sigma_t^{k}(\btheta_k,\bl_k)'\pi\biggr|^2 \\
&&\qquad\quad{} + \frac{1}{p} U(Y_t^k(\btheta_k,\bl_k)) \int_E U\bigl(\pi \gamma
_t^k(\btheta_k,\bl_k,\ell)\\
&&\qquad\quad\hspace*{112.6pt}{}
- Y_t^{k+1}(\btheta_k,t,\bl_k,\ell)\bigr) \eta_{{k+1}}(\bl_k,d\ell)
\end{eqnarray*}
for $k = 0,\ldots,n-1$ and $\hat g^n(\pi,t,\omega,\btheta,\bl) = |
Z^n_t(\btheta_k,\bl_k) + \frac{1}{p} \lambda_t^{n}(\btheta,\bl)
- \sigma_t^{n}(\btheta,\bl)'\pi|^2$. Recall that the indexed $\F$-adapted
processes $Y^k$ and $Y^{k+1}$ are continuous, hence $\F$-predictable.
Therefore, the map $(\pi,t,\omega,\btheta_k,\bl_k) \rightarrow\hat
g^k(\pi,t,\omega,\btheta_k,\bl_k)$ is $\Bc(\R^d)\otimes\Pc(\F
)\otimes\Bc(\Delta_k)\otimes\Bc(E^k)$-measurable.
Moreover, for $k = 0,\ldots,n$, $(\btheta_k,\bl_k) \in\Delta
_k\times E^k$, we recall from Remark~\ref{rempremium}
that either $\sigma^k(\btheta_k,\bl_k) = 0$ and $\gamma^k(\btheta
_k,\bl_k,\ell) = 0$, in which case, the continuous function
$\pi\rightarrow\hat g_k(\pi,t,\omega,\btheta_k,\bl_k)$ attains
trivially its infimum for $\pi= 0$, or $\sigma^k(\btheta_k,\btheta
_k)$ and
$\gamma^k(\btheta_k,\bl_k,\ell)$ are in the form $\sigma^k(\btheta_k,\bl
_k) = (\bar\sigma^k(\btheta_k,\bl_k) 0 )$,
$\gamma^k(\btheta_k,\bl_k,\ell) = (\bar\gamma^k(\btheta_k,\bl_k,\ell
) 0)$ for some full rank matrix $\bar\sigma^k(\btheta_k,\bl_k)$.
In this case, the infimum of $\hat g^k(\pi,\cdot)$ over $\pi\in A^k$
is equal to the infimum over $\bar\pi\in(\sigma^k)'A^k$ of
function $\bar g^k(\bar\pi,\cdot)$ where
\begin{eqnarray*}
\hspace*{-4pt}&&
\bar g^k(\bar\pi,t,\omega,\btheta_k,\bl_k) \\
\hspace*{-4pt}&&\qquad=
\frac{p}{2} \biggl| Z^k_t(\btheta_k,\bl_k) + \frac{1}{p} \lambda
_t^{k}(\btheta_k,\bl_k) - \bar\pi\biggr|^2 \\
\hspace*{-4pt}&&\qquad\quad{} + \frac{1}{p} U(Y_t^k) \int_E U\bigl( (\bar\sigma^k(\bar\sigma
^k)')^{-1}\bar\pi\bar\gamma_t^k(\ell)
- Y_t^{k+1}(\btheta_k,t,\bl_k,\ell)\bigr) \eta_{{k+1}}(\bl_k,d\ell)
\end{eqnarray*}
for $k = 0,\ldots,n-1$, and $\bar g^n(\hat\pi,t,\omega,\btheta,\bl
) = | Z^n_t(\btheta_k,\bl_k) + \frac{1}{p} \lambda_t^{n}(\btheta,\bl)
- \bar\pi|^2$. We clearly have
\[
\bar g_k(0,t,\omega,\btheta_k,\bl_k) \leq\liminf_{|\bar\pi|
\rightarrow\infty} \bar g_k(\bar\pi,t,\omega,\btheta_k,\bl_k),
\]
which shows that the continuous function
$\bar\pi\rightarrow\bar g^k(\bar\pi,t,\omega,\btheta_k,\bl_k)$
attains its infimum over the closed set
$(\sigma_t^k)'A^k$, and thus the function $\pi\rightarrow\hat
g^k(\pi,t,\omega,\btheta_k,\bl_k)$ attains its infimum over $A^k(\btheta
_k,\bl_k)$.
By a classical measurable selection theorem (see, e.g.,~\cite{wag80}),
one can then find for any $k = 0,\ldots,n$, $\hat\pi^k \in\Pc
_\F(\Delta_k,E^k)$ s.t.,
\[
\hat\pi_t^k(\btheta_k,\bl_k) \in\argmin_{\pi\in A^k(\btheta_k,\bl_k)}
\hat g^k(\pi,t,\btheta_k,\bl_k),\qquad \theta_k \leq t \leq T \qquad\mbox{a.s.}
\]
for all $(\btheta_k,\bl_k) \in\Delta_k(T)\times E^k$. Let us now
check that the trading strategy $\hat\pi= (\hat\pi^k)_{k=0,\ldots
,n}$ is admissible in the sense of Definition~\ref{defadmi}. First, by
writing that $\hat g^k(\hat\pi_t^k,t,\btheta_k,\bl_k) \leq\hat
g^k(0,t,\btheta_k,\bl_k)$, we get, similarly to \eqref{pileqz},
the existence of some constant $C$ s.t.
%
\begin{equation} \label{hatpiz}\qquad\quad
|\sigma_t^k(\btheta_k,\bl_k)'\hat\pi_t^k(\btheta_k,\bl_k)| \leq C\bigl( 1 +
|Z^k_t(\btheta_k,\bl_k)|\bigr),\qquad \theta_k\leq t\leq T \qquad\mbox{a.s.}
\end{equation}
for all $(\btheta_k,\bl_k) \in\Delta_k(T)\times E^k$, $k = 0,\ldots
,n$. Since $Z^k \in\Lb_W^2(\Delta_k,E^k)$,
and recalling~\ref{HB}(i), this shows that $\hat\pi^k$ satisfies \eqref
{integpi} for all $k = 0,\ldots,n$. Let us denote by $\hat X^{k,x}$
the wealth
process in \eqref{dynXI} controlled by $\hat\pi^k$, and starting from
$x$ at $\theta_k$.
By definition of $\hat\pi^k$, we have
%
\begin{eqnarray} \label{fkhatpi}
&&
f^{k}(t,Y^k_t,Z_t^k,\btheta_k,\bl_k) \nonumber\\
&&\qquad= \frac{p}{2} | Z_t^k -
\sigma_t^{k}(\btheta_k,\bl_k)'\hat\pi_t^k |^2
- b_t^{k}(\btheta_k,\bl_k)'\hat\pi_t^k \\
&&\qquad\quad{} + \frac{1}{p} U(Y_t^k) \int_E U\bigl(\hat\pi_t^k \gamma
_t^k(\btheta_k,\bl_k,\ell) - Y_t^{k+1}(\btheta_k,t,\bl_k,\ell)\bigr)
\eta_{{k+1}}(\bl_k,d\ell)\hspace*{-18pt}  \nonumber
\end{eqnarray}
for $k = 0,\ldots,n-1$, and $f^{n}(t,Y^n_t,Z_t^n,\btheta,\bl) = \frac
{p}{2} | Z_t^n - \sigma_t^{n}(\btheta,\bl)'\hat\pi_t^n
|^2 - b_t^{n}(\btheta,\bl)'\hat\pi_t^n$ for $k = n$. From the
forward dynamics of $Y^k$, we can then write for all $\theta_k\leq
t\leq T$
\[
U(\hat X_{t}^{k,x} - Y^k_{t}) = U(x -Y_{\theta_k}^k)
\Ec_{t}^k \bigl( p\bigl(Z^k- (\sigma^k)'\hat\pi^k\bigr)\bigr)R_{t}^k
\]
with
\begin{eqnarray*}
\Ec_{t}^k \bigl( p\bigl(Z^k\!-\! (\sigma^k)'\hat\pi^k\bigr)\bigr)\! =\!
\exp\biggl(\! p \!\int_{\theta_k}^t\! \bigl(Z^k_{s}\!-\! (\sigma^k_{s})'\hat\pi_s^k\bigr)\,dW_{s}
\!-\!\frac{p^2}{2}\! \int_{\theta_k}^t\! |Z^k_{s}- (\sigma^k_{s})'\hat\pi
_s^k|^2 \,ds\! \biggr)
\end{eqnarray*}
and
\begin{eqnarray*}
R_{t}^k &=& \exp\biggl(-\int_{\theta_k}^t U(Y_{s}^k)
\int_E U\bigl(\hat\pi_t^k \gamma_t^k(\btheta_k,\bl_k,\ell)\\
&&\hspace*{104pt}{} -
Y_t^{k+1}(\btheta_k,t,\bl_k,\ell)\bigr)
\eta_{{k+1}}(\bl_k,d\ell) \,ds\biggr)
\end{eqnarray*}
for $k = 0,\ldots,n-1$, and $R_t^n = 1$. Now, from step 2 and
\eqref{hatpiz}, the process
$\int_{\theta_k}^{ \cdot} p(Z^k - (\sigma^k)'\hat\pi^k) \,dW$ is a
BMO-martingale, and hence (see~\cite{kaz}),
$\Ec^k(p(Z^k - (\sigma^k)'\hat\pi^k))$ is of class (D).
Moreover, since $U$ is nonpositive, we see that
$|R^k| \leq1$, and so $|U(\hat X^{k,x}-Y^k)| \leq U(x-Y_{\theta
_k}^k) \Ec^k(p(Z^k - (\sigma^k)'\hat\pi^k))$,
which shows that $U(\hat X^{k,x}-Y^k)$ is of class (D), and then also
$U(\hat X^{k,x})$ since $Y^k$ is essentially bounded.
It remains to check that for all $k = 0,\ldots,n-1$, \mbox{$(\btheta_k,\bl
_k) \in\Delta_k(T)\times E^k$},
\[
\mathbb{E}\biggl[\int_{\theta_k}^T \int_{E} (-U)\bigl(\hat X^{k,x}_{t}+
\hat\pi^k_{t} \gamma^k_{t}(\btheta_{k},\bl_{k},\ell) \bigr)
\eta_{{k+1}}(\bl_k,d\ell) \,ds\biggr] < \infty.
\]
By the definition of $\hat\pi^k$ [which implies \eqref{fkhatpi}], the
process $\xi^k(x,\btheta_k,\bl_k,\hat\pi^k)$ defined in step 1, is a
local martingale.
By considering a localizing $\F$-stopping times sequence $(\rho_n)_n$
valued in $[\theta_k,T]$ for this local martingale, we obtain
\begin{eqnarray*}
& & \Ee\biggl[ \int_{\theta_k}^{T\wedge\rho_n}
\int_E (-U)\bigl(\hat X_t^{k,x} + \hat\pi_t^k \gamma_t^k(\ell) -
Y_t^{k+1}(\btheta_k,t,\bl_k,\ell) \bigr) \eta_{{k+1}}(\bl_k,d\ell)\,dt
\biggr] \\
&&\qquad= \Ee[ U(\hat X_{T\wedge\rho_n}^{k,x} - Y_{T\wedge\rho_n}^k) -
U(x-Y_{\theta_k}^k) ] \leq
\Ee[ -U(x-Y_{\theta_k}^k) ],
\end{eqnarray*}
since $U$ is nonpositive. By Fatou's lemma, we get the required
inequality, and this proves that $\hat\pi^k \in\Ac_\F^k$, for any
$k = 0,\ldots,n$; that is, $\hat\pi= (\hat\pi^k)_{k=0,\ldots,n}$ is
admissible: $\hat\pi\in\Ac_\G$.\vspace*{10pt}

\textit{Step} 4: Since $\hat\pi= (\hat\pi^k)_{k=0,\ldots
,n}$ is admissible, and recalling that the processes $Y^k$ are
essentially bounded, this implies that the local martingales $\xi
^k(x,\btheta_k,\bl_k,\hat\pi^k)$, $k = 0,\ldots,n$, are ``true''
martingales. Hence, the inequalities
in \eqref{Unsuper}--\eqref{Uksuper} become equalities for $\nu= \hat
\pi$, which yield
%
\begin{eqnarray}
\label{Unmar}
U\bigl(\hat X_t^{n,x} - Y_{t}^{n}(\btheta,\bl)\bigr) & = &
\Ee\bigl[ U\bigl(\hat X_T^{n,x} - H_T^n(\btheta,\bl)\bigr) \alpha_T(\btheta,\bl)
| \Fc_t \bigr],\\[-25pt]\nonumber
\end{eqnarray}
\begin{eqnarray}
\label{Ukmar}
&&
U\bigl(\hat X_t^{k,x} - Y_{t}^{k}(\btheta_k,\bl_k)\bigr)
\nonumber\\
&&\qquad =
\Ee\biggl[ U\bigl(\hat X_T^{k,x} - H_T^k(\btheta_k,\bl_k)\bigr) \alpha_T^k(\btheta
_k,\bl_k) \nonumber\\[-8pt]\\[-8pt]
&&\hspace*{11.6pt}\qquad\quad{} + \int_{t}^T \int_E U\bigl(\hat X_s^{k,x} + \hat\pi_s^k \gamma
_s^k(\ell)\nonumber\\
&&\hspace*{67.5pt}\qquad\quad{}  - Y_s^{k+1}(\btheta_k,s,\bl_k,\ell) \bigr) \eta_{{k+1}}(\bl
_k,d\ell) \,ds
\Big| \Fc_t \biggr] \nonumber
\end{eqnarray}
for $k = 0,\ldots,n$, $(\btheta_k,\bl_k) \in\Delta_k(T)\times
E^k$, $t \in[\theta_k,T]$, $x \in\R$. Let us prove the properties
\eqref{VkYk} by backward induction on $k = 0,\ldots,n$.
For $k = n$, from the additive form of the wealth process $X^{n,x}$
and the exponential form of the utility function~$U$,\vadjust{\goodbreak} we observe that for
any $t \in[\theta_n,T]$, $\pi^n \in\Ac_\F(t,\nu^n)$, the quantity
\[
\Ee\biggl[ \frac{U(X_T^{n,x} - H_T^n(\btheta,\bl))}{-U(X_t^{n,x})} \alpha
_T(\btheta,\bl) \Big| \Fc_t \biggr]
\]
does not depend on the choice $\nu^n \in\Ac_\F^n$. By combining
\eqref{Unsuper} and \eqref{Unmar}, we then have
\begin{eqnarray*}
J_t^n(\btheta,\bl) :\!&= & \esssup_{\pi^n\in\Ac_\F^n(t,\nu^n)} \Ee\biggl[
\frac{U(X_T^{n,x} - H_T^n(\btheta,\bl))}{-U(X_t^{n,x})} \alpha_T(\btheta
,\bl) \Big| \Fc_t \biggr] \\
& \leq& U(- Y_{t}^{n}(\btheta,\bl)) =
\Ee\biggl[ \frac{U(\hat X_T^{n,x} - H_T^n(\btheta,\bl))}{-U(\hat
X_t^{n,x})} \alpha_T(\btheta,\bl) \Big| \Fc_t \biggr] \\
&\leq&
J_t^n(\btheta,\bl),
\end{eqnarray*}
where we used in the last inequality the trivial fact that $\hat\pi^n
\in\Ac_\F^n(t,\hat\pi^n)$. This shows that
$U(- Y_{t}^{n}(\btheta,\bl)) = J_t^n(\btheta,\bl)$, and so
$V^n_t(x,\btheta,\bl,\nu^n) = U(X_t^{n,x}-Y_t^n(\btheta,\bl))$ for
any $\nu^n \in\Ac_\F^n$, $x \in\R$, $(\btheta,\bl) \in\Delta
_n(T)\times E^n$, which is property \eqref{VkYk} at step $k = n$.
Assume now that \eqref{VkYk} holds true at step $k+1$. Then, we observe,
similarly as above, that
for any $t \in[\theta_k,T]$, $\pi^k \in\Ac_\F(t,\nu^k)$, the quantity
\begin{eqnarray*}
& & \Ee\biggl[ \frac{U(X_T^{k,x} - H_T^{k}(\btheta_k,\bl_k) )
\alpha_T^k(\btheta_k,\bl_k)}{-U(X_t^{k,x})} \\[-0.5pt]
& &\quad\hspace*{0pt}{} + \int_{t}^T \int_E \frac{V^{k+1}_{\theta_{k+1}}
( X_{\theta_{k+1}}^{k,x} + \pi_{\theta_{k+1}}^k \gamma_{\theta
_{k+1}}^k(\ell_{k+1}),
\btheta_{k+1},\bl_{k+1})}{-U(X_t^{k,x})} \\[-0.5pt]
&&\qquad\quad\hspace*{95.3pt}{}\times\eta_{{k+1}}(\bl_k,d\ell
_{{k+1}}) \,d\theta_{k+1} \Big| \Fc_{t} \biggr] \\[-0.5pt]
&&\qquad= \Ee\biggl[ \frac{U(X_T^{k,x} - H_T^{k}(\btheta_k,\bl_k) )
\alpha_T^k(\btheta_k,\bl_k)}{-U(X_t^{k,x})} \\[-0.5pt]
&&\qquad\qquad{} + \int_{t}^T \int_E \frac{U( X_{s}^{k,x} + \pi
_{s}^k \gamma_{s}^k(\ell) -Y_s^{k+1}(\btheta_k,s,\bl_l,\ell)
)}{-U(X_t^{k,x})}\\[-0.5pt]
&&\qquad\quad\hspace*{141pt}{}\times
\eta_{{k+1}}(\bl_k,d\ell) \,ds \Big| \Fc_{t} \biggr]
\end{eqnarray*}
is independent of the choice $\nu^k \in\Ac_\F^k$. By combining
\eqref{Uksuper} and \eqref{Ukmar}, we then have
\begin{eqnarray*}
&&
J_t^k(\btheta_k,\bl_k) \\[-2pt]
&&\qquad:= \esssup_{\pi^k\in\Ac_\F^k(t,\nu^k)}
\Ee\biggl[ \frac{U(X_T^{k,x} - H_T^{k}(\btheta_k,\bl_k) ) \alpha
_T^k(\btheta_k,\bl_k)}{-U(X_t^{k,x})} \\[-2pt]
&&\hspace*{62pt}\qquad\quad{}+ \int_{t}^T \int_E \frac{V^{k+1}_{\theta_{k+1}} ( X_{\theta
_{k+1}}^{k,x} + \pi_{\theta_{k+1}}^k \gamma_{\theta_{k+1}}^k(\ell_{k+1}),
\btheta_{k+1},\bl_{k+1})}{-U(X_t^{k,x})} \\[-2pt]
&&\qquad\quad\hspace*{179.6pt}{}\times\eta_{{k+1}}(\bl_k,d\ell
_{k+1}) \,d\theta_{k+1} \Big| \Fc_{t} \biggr] \\[-4pt]
&&\qquad \leq U(- Y_{t}^{k}(\btheta_k,\bl_k)) \\[-2pt]
&&\qquad= \Ee\biggl[ \frac{U(\hat X_T^{k,x} - H_T^{k}(\btheta_k,\bl_k)
) \alpha_T^k(\btheta_k,\bl_k)}{-U(X_t^{k,x})} \\[-2pt]
&&\qquad\quad\hspace*{11.2pt}{}+ \int_{t}^T \int_E \frac{V^{k+1}_{\theta_{k+1}} ( \hat X_{\theta
_{k+1}}^{k,x} + \hat\pi_{\theta_{k+1}}^k \gamma_{\theta
_{k+1}}^k(\ell_{k+1}),
\btheta_{k+1},\bl_{k+1})}{-U(\hat X_t^{k,x})}\\[-2pt]
&&\qquad\quad\hspace*{128.6pt}{}\times \eta_{{k+1}}(\bl
_k,d\ell_{k+1}) \,d\theta_{k+1} \Big| \Fc_{t} \biggr] \\[-4pt]
&&\qquad \leq J_t^k(\btheta_k,\bl_k),
\end{eqnarray*}
where we used in the last inequality the trivial fact that $\hat\pi^k
\in\Ac_\F^k(t,\hat\pi^k)$. This proves that
$U(- Y_{t}^{k}(\btheta_k,\bl_k)) = J_t^k(\btheta_k,\bl_k)$,
and thus the property \eqref{VkYk} at step $k$. Notice that this representation
of $Y^k$ shows as a byproduct the uniqueness of the solution to the
recursive system of BSDEs (\ref{EK}). Finally, relations
\eqref{Ukmar} for $t = \theta_k$, together with \eqref{VkYk}, yield
\begin{eqnarray*}
V^{n}(x,\btheta,\bl) &=& \Ee[ U(\hat X_T^{n,x} - H_T^{n})
\alpha_T(\btheta,\bl) | \Fc_{\theta_{n}} ], \\[-2pt]
V^k(x,\btheta_k,\bl_k) &=& \Ee\biggl[ U(\hat X_T^{k,x} - H_T^{k}
) \alpha_T^k(\btheta_k,\bl_k)\\[-2pt]
& &\hspace*{10.7pt}{}+ \int_{\theta_{k}}^T \int_E V^{k+1}\bigl( \hat X_{\theta_{k+1}}^{k,x}
+ \hat\pi_{\theta_{k+1}}^k \gamma_{\theta_{k+1}}^k(\ell_{k+1}),
\btheta_{k+1},\bl_{k+1}\bigr)\\[-2pt]
&&\hspace*{121pt}{}\times \eta_{{k+1}}(\bl_k,d\ell_{k+1}) \,d\theta
_{k+1} \Big| \Fc_{\theta_{k}} \biggr],
\end{eqnarray*}
which prove that $\hat\pi= (\hat\pi^k)_{k=0,\ldots,n}$ is an
optimal trading strategy.
\end{pf}
%
\begin{Rem}
We recall that, in a default-free market, the It\^o model for
stock price $S$ with risk premium $\lambda$ and volatility $\sigma$,
the optimal trading strategy (in amount) for an exponential utility
function $U(x) = -e^{-px}$, and option payoff $H_T$, is given by
(see~\cite{HIM} or~\cite{REK})
\[
\hat\pi_t^M \in\argmin_{\pi\in A}
\biggl| Z_t + \frac{1}{p} \lambda_t - (\sigma_t)'\pi\biggr|^2,
\]
where $(Y,Z)$ is the solution to the BSDE $dY_t = -f(t,Z_t)\,dt +
Z_TdW_t$, $Y_T = H_T$, $f(t,z) = \inf_{\pi\in A}|z+\frac{1}{p}\lambda_t
- (\sigma_t)'\pi|^2$. In\vspace*{1pt} our multiple defaults risk model,
inducing jumps on the stock price, we see from \eqref{formpik} the
influence of jumps in the optimal\vadjust{\goodbreak} trading strategy $\hat\pi^k$ within
the $k$-default scenario: there is a similar term involving the
coefficients $\lambda^k$ and $\sigma^k$ corresponding to the
default-free regime case, but the investor will take into account the
possibility of a default and jump before the final horizon, and which
is formalized by the additional term involving the jump
size~$\gamma^k$. In particular, if $\gamma^k$ is negative (in the
one-asset case $d = 1$), meaning that there is a loss at default. Then
the infimum in \eqref{formpik} will be achieved for a value $\hat\pi^k$
smaller than the one without jumps. This means that when the investor
knows that there will be a loss at default on the stock, he will invest
less in this asset, which is intuitive. In the next section, we shall
measure quantitatively this impact on a two-assets model with 
defaults.\vspace*{-2pt}
\end{Rem}

\section{Applications and numerical illustrations}\label{sec5}

For numerical illustrations, we consider a portfolio of two defaultable names,
and denote by $\tau_1$ and $\tau_2$ their respective nonordered default times,
assumed to be independent of $\bF$, so that their conditional density
(w.r.t. $\bF$) is a deterministic function.
We suppose that $\tau_1$ and $\tau_2$ are correlated via the Gumbel
copula which is suitable to characterize heavy tail dependence and is
often used for insurance portfolios. More precisely, we let $\proba[\tau
_1>\theta_1,\tau_2>\theta_2|\cF_t] =
\proba[\tau_1>\theta_1,\tau_2>\theta_2] = \exp(-((a_1\theta
_1)^\beta+(a_2\theta_2)^\beta)^{1/\beta})$ with
$a_1,a_2>0$ and $\beta\geq1$. In this model, each marginal default
time $\tau_i$ follows the exponential law with constant intensity
$a_i$, $i = 1,2$, and the correlation between the two defaults is
characterized by the constant parameter $\beta$. The case $\beta=1$
corresponds to the independence case, and a larger value of $\beta$
implies a large linear correlation between the survival events
$\rho^s(T) = \operatorname{corr}(\indic_{\{\tau_1>T\}},\indic_{\{\tau
_2> T\}})$.
The default density of
$\btau= (\tau_1,\tau_2)$ is thus given by
\[
\alpha^{\btau}(\theta_1,\theta_2)
=G(\theta_1,\theta_2)\frac{(a_1a_2)^\beta}{(\theta_1\theta_2)^{1-\beta
}}
u(\theta_1,\theta_2)^{1-2\beta}\bigl(u(\theta_1,\theta_2)+\beta-1\bigr),
\]
where $G(\theta_1,\theta_2)=\proba(\tau_1>\theta_1,\tau_2>\theta_2)=\exp
(-u(\theta_1,\theta_2))$. As explained in Section~\ref{sec2.1}
and Remark~\ref{remdens},
the case of ordered default times $\hat\tau_1 = \min(\tau_1,\tau
_2)$, $\hat\tau_2 = \max(\tau_1,\tau_2)$ can be recovered by
considering the marks $(\iota_1,\iota_2)$ indicating the order of the
defaults $(\tau_1,\tau_2)$. The density
of $(\hat\tau_1,\hat\tau_2,\iota_1,\iota_2)$ is given by
\[
\alpha(\btheta,i,j) = \indic_{\{i=1,j=2\}} \alpha^{\btau}(\theta
_1,\theta_2) + \indic_{\{i=2,j=1\}}
\alpha^{\btau}(\theta_2,\theta_1)
\]
for $\btheta= (\theta_1,\theta_2) \in\Delta_2$. Before any
default, the price process $S^0 = (S^{1,0},S^{2,0})$ of the two
names is governed by a two-dimensional Black--Scholes model with the correlation
\[
dS_t^0 = S_t^0* (b^0 \,dt + \sigma^0 \,dW_t),
\]
where $b^0 = (b^{1,0},b^{2,0})$ is a constant vector in $\R^2$,
$\sigma^0$ is the constant matrix
\[
\sigma^0 =
\pmatrix{
\sigma^{1,0}\sqrt{1-\rho^2} & \sigma^{1,0}\rho\cr
0 & \sigma^{2,0}}
\]
with $\sigma^{1,0} > 0$, $\sigma^{2,0} > 0$, $\rho\in(-1,1)$ and $W =
(W^1,W^2)$
is a two-dimensional Brownian motion. The associated risk premium\vadjust{\goodbreak} is
then given by $\lambda^0 = (\lambda^{1,0},\break\lambda^{2,0})$ with
\[
\lambda^{1,0} = \frac{1}{\sqrt{1-\rho^2}} \biggl( \frac
{b^{1,0}}{\sigma^{1,0}} - \rho\frac{b^{2,0}}{\sigma^{2,0}}
\biggr),\qquad
\lambda^{2,0} = \frac{b^{2,0}}{\sigma^{2,0}}.
\]
Once the name $j$ defaults at time $\tau_j$, it drops to zero, but it
also incurs a constant relative jump (loss\vspace*{1pt} or gain) of
size $\gamma^i \in[-1,\infty)$ on the other name $i \neq j$. We denote
by $S^{i,1}(\theta_1) = S^{i,1}(\theta_1,j)$ the price process of the
survival name $i$ after the first default due to name $j \neq i$ at
time $\tau_j = \theta _1$. We then have $S^{i,1}_{\theta_1}(\theta_1) =
S^{i,0}_{\theta _1}(1+\gamma^i)$, and we assume that it follows a
Black--Scholes model
\[
dS_t^{i,1}(\theta_1) = S_t^{i,1}(\theta_1)(b^{i,1} \,dt + \sigma^{i,1}
\,dB_t^i),\qquad t \geq\theta_1,
\]
with\vspace*{1pt} constants $b^{i,1}$ and $\sigma^{i,1} > 0$. Here
$B^i$ is the Brownian motion $B^1 = \sqrt{1-\rho^2}W^1+\rho W^2$, $B^2
= W^2$. Finally, after both defaults, the two names cannot be traded
anymore, that is, $S^2 = (S^{1,2},S^{2,2}) = 0$.

We consider the investment problem with utility function $U(x) =
-e^{-px}$, without option payoff $H_T = 0$,
without portfolio constraint, and solve the recursive system of BSDEs.
Since all the coefficients of the assets price and the density are
deterministic, we notice that these BSDEs reduce actually to ordinary
differential equations (ODEs).
We start from the case $n = 2$ after the defaults of both names.
The solution to the BSDE (\ref{En})
for $n = 2$ is clearly degenerate:
\[
Y^2(\btheta,i,j) = \frac{1}{p} \ln\alpha(\btheta,i,j),\qquad \btheta
= (\theta_1,\theta_2) \in\Delta_2, i,j \in\{1,2\}, i\neq j.
\]
Let us denote by $Y^{1,i}(\theta_1) = Y^1(\theta_1,i)$, $i = 1,2$, the
solution to the BSDE (E1) after the first default
due to
name $i$. Notice that the {auxiliary function $\alpha^{1,i}(\theta_1) =
\alpha^1(\theta_1,i)$, defined in \eqref{AlphaAuxk}}, is given for
$ i,j=1,2,
i\neq j$, by
\begin{eqnarray*}
\alpha^{1,i}_t(\theta_1) &=& \int_t^\infty\alpha(\theta_1,\theta_2,i,j)
\,d\theta_2 \\
&=& \frac{a_i^\beta}{\theta_1^{1-\beta}}
\bigl((a_i\theta_1)^\beta+(a_jt)^\beta\bigr)^{1/\beta}
e^{-((a_i\theta_1)^\beta+(a_jt)^\beta)^{1/\beta}}.
\end{eqnarray*}
The function $Y^{1,i}$ is then given by the solution to the ODE
\begin{eqnarray*}
Y_t^{1,i}(\theta_1) &=&
\frac{1}{p}\biggl[\beta\ln a_i+(\beta-1)\ln\theta_1\\
&&\hspace*{12.9pt}{}+\frac{1}{\beta}\ln
\bigl((a_i\theta_1)^\beta+(a_jt)^\beta\bigr)
-\bigl((a_i\theta_1)^\beta+(a_jt)^\beta\bigr)^{1/\beta} \biggr]\\
&&{} + \int_t^T f^{1,i}(s,Y_s^{1,i},\theta_1) \,ds,
\end{eqnarray*}
where
\begin{eqnarray*}
f^{1,i}(t,y,\theta_1) &=& - \frac{1}{2p}\biggl|\frac{b^{j,1}}{\sigma
^{j,1}}\biggr|^2
+ \inf_{\pi\in\R} \biggl\{ \frac{p}{2} \biggl| \frac{1}{p} \frac
{b^{j,1}}{\sigma^{j,1}} - \sigma^{j,1}\pi\biggr|^2
+ \frac{1}{p} e^{-p(y-\pi)} \alpha(\theta_1,t,i,j) \biggr\}
\end{eqnarray*}
for $i,j \in\{1,2\}$, $i \neq j$. For $k = 0$, the
survival probability $\alpha^0$ is equal to
\[
\alpha_T^0 = \P[ \tau_1 >T,\tau_2 > T] = \exp\bigl(- T(a_1^\beta
+a_2^\beta)^{1/\beta} \bigr),
\]
and the function $Y^0$ to the BSDE (E0) is then given by the
solution to the ODE
%
\begin{equation}\label{Yt0}
Y_t^{0} = - \frac{T}{p}(a_1^\beta+a_2^\beta)^{1/\beta} + \int_t^T
f^{0}(s,Y_s^{0}) \,ds,
\end{equation}
where
\begin{eqnarray*}
f^0(t,y) & = & -\frac{1}{2p}|\lambda^0|^2\\
&&{} + \inf_{\pi=(\pi^1,\pi^2) \in
\R^2} \biggl\{ \frac{p}{2} \biggl|\frac{1}{p} \lambda^0 - (\sigma^0)'\pi
\biggr|^2 \\
& &\hspace*{73.2pt}{} + \frac{1}{p} e^{-py} \bigl[ e^{-p(-\pi^1 +
\pi^2\gamma^2 - Y_t^{1,1}(t))} \\
&&\hspace*{119pt}{}+
e^{-p(\pi^1\gamma^1 - \pi^2 - Y_t^{1,2}(t))}
\bigr] \biggr\}.
\end{eqnarray*}

We perform numerical results to study notably the following parameters:
the loss or gain at default, the default intensities and the
correlation between the defaults and between the assets. We choose the
parameters of assets as below and fix them to be the same in all our
tests: {$b^{1,0}=b^{2,0}=0.02$, $\sigma^{1,0}=\sigma^{2,0}=0.1$,
$b^{1,1}=b^{2,1}=0.01$, $\sigma^{1,1}=\sigma^{2,1}=0.2$,} $p=1$ and $T=1$.

%
\begin{figure}

\includegraphics{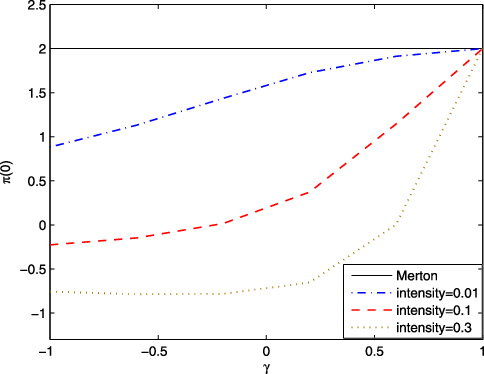}\vspace*{-2pt}

\caption{Optimal strategy $\hat\pi$ before any default vs Merton $\hat\pi^{M}$.}
\label{intensitypi}\vspace*{-3pt}
\end{figure}

\begin{figure}[b]
\vspace*{-3pt}
\includegraphics{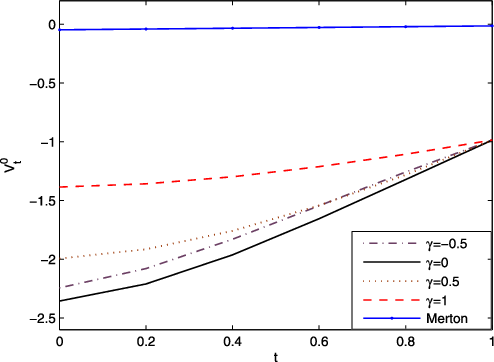}\vspace*{-2pt}

\caption{Value function $V_t^0$.}
\label{gammaYt}
\end{figure}

In Figure~\ref{intensitypi}, we present the optimal strategies $\hat\pi
=(\hat\pi^1,\hat\pi^2)$ at the initial time before any default, for
different values of loss or gain at default and of default intensity.
In Figure~\ref{intensitypi}, we consider a symmetric case where the
default intensities $a_1$ and $a_2$, and the loss/gain $\gamma^1$ and
$\gamma^2$, are equal, respectively, so they are the same for $\hat\pi
^1$ and $\hat\pi^2$. We choose the correlation parameter $\rho=0$ and
$\beta=2$. The optimal strategy is increasing with respect to $\gamma$,
which means that one should invest less on the assets when there is a
large loss of default. When $\gamma=1$, the strategy converges to the
Merton one, since in this case, the gain at default of the surviving
name will recompense the total loss of the default one. Furthermore,
the strategy is decreasing with respect to the default intensity. So
when there is a higher risk of default, one should reduce her
investment. In particular, if the default probability is high, and the
loss at default is large, then the investor should sell instead of buy
the assets. Only when $\gamma$ becomes positive, and the gain at
default is large enough to recompense the default risks, she can choose
to buy the asset again.\vadjust{\goodbreak}

%
\begin{figure}

\includegraphics{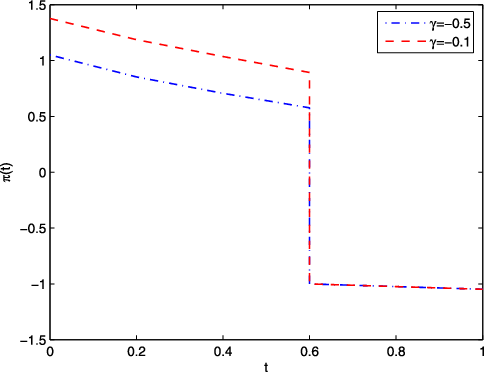}

\caption{Time evolution of the optimal strategy $\hat\pi$ given a default.}
\label{pit}\vspace*{-3pt}
\end{figure}

%
\begin{table}[b]\vspace*{-3pt}
\def\arraystrech{0.9}
\caption{Optimal strategies ${\hat\pi}^{1}$ and ${\hat
\pi}^{2}$ before any defaults with various $\gamma$ and default
intensities}\label{Tabintensitypi}
\begin{tabular*}{\tablewidth}{@{\extracolsep{\fill}}ld{2.3}d{2.3}d{2.3}d{2.3}
d{1.3}c@{}}
\hline
& \multicolumn{5}{c}{$\bolds{\gamma}$}\\[-4pt]
& \multicolumn{5}{c}{\hrulefill}\\
&\multicolumn{1}{c}{$\bolds{-0.5}$}&\multicolumn{1}{c}{$\bolds{-0.1}$}
& \multicolumn{1}{c}{$\bolds{0}$} &\multicolumn{1}{c}{$\bolds{0.5}$}
& \multicolumn{1}{c}{$\bolds{1}$} & \multicolumn{1}{c@{}}{\textbf{Merton}}\\
\hline
$a_1=0.01$, $a_2=0.1$, $\beta=2$&\multicolumn{1}{c}{$\rho^s=0.2936$}\\
[3pt]
\quad${\hat\pi}^1$ & 0.462 & 1.659 & 1.892 & 2.621 & 2.832 & 2\\
\quad${\hat\pi}^2$ & -1.047 & -0.709 & -0.498 & 0.623 & 1.168
&2\\
[3pt]
$a_1=0.1$, $a_2=0.1$, $\beta=2$&\multicolumn{1}{c}{$\rho^s=0.5736$}\\
[3pt]
\quad${\hat\pi}^1$ & -0.353 & -0.210 & -0.147 & 0.556 & 2 & 2\\
\quad${\hat\pi}^2$ & -0.353 & -0.210 & -0.147 & 0.556 & 2 & 2\\
[3pt]
$a_1=0.3$, $a_2=0.1$, $\beta=2$& \multicolumn{1}{c}{$\rho^s=0.4555$}\\
[3pt]
\quad${\hat\pi}^1$ & -1.723 & -1.719 & -1.647 & -0.697 & 1.293 &
2\\
\quad${\hat\pi}^2$ & -0.132 & 0.453 & 0.521 & 1.121 & 2.707 &2\\
\hline
\end{tabular*}
\end{table}

Figure~\ref{gammaYt} plots the evolution of the value function before
default, that is, $t \rightarrow V_t^0(x)=-e^{-p(x-Y_t^0)}$, where
$Y_t^0$ is the solution of equation \eqref{Yt0}, and we have chosen
$x=0$ in the test. We consider various values of $\gamma$ with the same
parameters as above and let $a_1=a_2=0.01$, $\beta=2$. The survival
correlation is equal to $\rho^s(T)=0.5846$. We observe a larger value
function when the gain at default
($\gamma>0$) is larger. We also notice that the value function in a
loss at default ($\gamma<0$) situation outperforms the no-loss
case ($\gamma= 0$), which means that one can take profit from a
loss of the risky stock by a shortsale strategy.

Figure~\ref{pit} plots the evolution of the optimal investment strategy
$\hat\pi(t)$ for $t \in[0,T]$, $T = 1$, when there is a
default event at time $\tau=0.6$, the parameters being the same as in
Figure~\ref{gammaYt}, with two different levels of loss at default
$\gamma$. We observe a jump of the trading strategy at the default time
in both curves. When there is a larger loss at default, one should
invest less from the beginning; however, after the default occurs, the
trading strategies on the surviving firm become the same whatever the
loss at default is.

We present, in Table~\ref{Tabintensitypi}, the optimal strategies at
initial time before defaults for firms with different levels of default
risks ($a_1\neq a_2$). We still suppose equal loss or gain at default
($\gamma^1=\gamma^2$). Similarly to Figure~\ref{intensitypi}, when the
default intensity $a_1$ of the first firm increases, one should reduce
the investment on this firm. In the case of high default risks and loss
at default, one should sell instead of buy the risky asset.
However, the strategy on the second firm (the one with $a_2=0.1$) will
in general increase when its counterparty becomes more
risky.

{Finally, we examine the impact of correlation parameters $\rho$ and
$\beta$ on the trading strategies before any default. In the following
test presented in Table~\ref{Tabcorrelpi}, we fix $a_1=0.01$ and
$a_2=0.1$. We observe that the correlation $\rho$ between the assets
will modify the benchmark Merton strategies. When $\rho$ increases, the
investment on the less risky asset goes in two directions: one should
increase its quantity in the loss at default case and reduce it in the
gain at default case; as for the more risky asset, one should always
reduce the investment. Concerning the parameter $\beta$, when there is
a larger $\beta$ and hence a higher correlation between the survival
events, one should increase the investment in the less risky asset and
decrease the investment in the more risky one.\vspace*{-3pt}
%
\begin{table}
\def\arraystrech{0.9}
\caption{Optimal strategies ${\hat\pi}^{1}$ and ${\hat\pi
}^{2}$ with various $\rho$ and $\beta$}\label{Tabcorrelpi}
\begin{tabular*}{\tablewidth}{@{\extracolsep{\fill}}ld{2.3}d{2.3}d{2.3}d{2.3}
d{1.3}d{1.3}@{}}
\hline
& \multicolumn{5}{c}{$\bolds{\gamma}$}\\[-4pt]
& \multicolumn{5}{c}{\hrulefill}\\
&\multicolumn{1}{c}{$\bolds{-0.5}$}&\multicolumn{1}{c}{$\bolds{-0.1}$}
& \multicolumn{1}{c}{$\bolds{0}$} &\multicolumn{1}{c}{$\bolds{0.5}$}
& \multicolumn{1}{c}{$\bolds{1}$} & \multicolumn{1}{c@{}}{\textbf{Merton}}\\
\hline
$\rho=0$, $\beta=1$&\multicolumn{1}{c}{$\rho^s=0$}\\
[3pt]
\quad${\hat\pi}^1$ & 0.228 & 0.942 & 1.099 & 1.966 & 2.459 & 2\\
\quad${\hat\pi}^2$ & -0.867 & -0.452 & -0.278 & 0.856 & 1.541 &
2\\
[3pt]
$\rho=0$, $\beta=2$&\multicolumn{1}{c}{$\rho^s=0.2936$}\\
[3pt]
\quad${\hat\pi}^1$ & 0.462 & 1.659 & 1.892 & 2.621 & 2.832 & 2\\
\quad${\hat\pi}^2$ & -1.047 & -0.709 & -0.498 & 0.623 & 1.168
&2\\
[3pt]
$\rho=0.3$, $\beta=1$&\multicolumn{1}{c}{$\rho^s=0$}\\
[3pt]
\quad${\hat\pi}^1$ & 0.492 & 1.081 & 1.188 & 1.715 & 2.025 &
1.539\\
\quad${\hat\pi}^2$ & -0.959 & -0.504 & -0.348 & 0.519 & 1.052
&1.539\\
[3pt]
$\rho=0.3$, $\beta=2$&\multicolumn{1}{c}{$\rho^s=0.2936$}\\
[3pt]
\quad${\hat\pi}^1$ & 0.863 & 1.939 & 2.077 & 2.399 & 2.450 &
1.539\\
\quad${\hat\pi}^2$ & -1.235 & -0.817 & -0.626 & 0.216 & 0.627
&1.539\\
\hline
\end{tabular*}     \vspace*{-3pt}
\end{table}



\printaddresses


\begin{thebibliography}{17}

\bibitem{ABE}
\begin{barticle}[mr]
\bauthor{\bsnm{Ankirchner},~\bfnm{Stefan}\binits{S.}},
  \bauthor{\bsnm{Blanchet-Scalliet},~\bfnm{Christophette}\binits{C.}} \AND
  \bauthor{\bsnm{Eyraud-Loisel},~\bfnm{Anne}\binits{A.}}
(\byear{2010}).
\btitle{Credit risk premia and quadratic {BSDE}s with a single jump}.
\bjournal{Int. J. Theor. Appl. Finance}
\bvolume{13}
\bpages{1103--1129}.
\bid{doi={10.1142/S0219024910006133}, issn={0219-0249}, mr={2738764}}
\bptok{imsref}%
\end{barticle}
\endbibitem

\bibitem{BC2008}
\begin{barticle}[mr]
\bauthor{\bsnm{Brigo},~\bfnm{Damiano}\binits{D.}} \AND
  \bauthor{\bsnm{Chourdakis},~\bfnm{Kyriakos}\binits{K.}}
(\byear{2009}).
\btitle{Counterparty risk for credit default swaps: Impact of spread volatility
  and default correlation}.
\bjournal{Int. J. Theor. Appl. Finance}
\bvolume{12}
\bpages{1007--1026}.
\bid{doi={10.1142/S0219024909005567}, issn={0219-0249}, mr={2574492}}
\bptok{imsref}%
\end{barticle}
\endbibitem

\bibitem{CJZ}
\begin{bmisc}[mr]
\bauthor{\bsnm{Cr{\' e}pey},~\bfnm{S.}\binits{S.}},
  \bauthor{\bsnm{Jeanblanc},~\bfnm{M.}\binits{M.}} \AND
  \bauthor{\bsnm{Zargari},~\bfnm{B.}\binits{B.}}
(\byear{2010}).
\bhowpublished{Counterparty risk on a CDS in a Markov chain copula model with joint defaults.
In \textit{Recent Advances in Financial Engineering} (M.~Kijima, C.~Hara, Y.~Muromachi and K.~Tanaka, eds.) 91--126. World Scientific, Singapore.}
\bptok{imsref}%
\end{bmisc}
\endbibitem

\bibitem{ejj1}
\begin{barticle}[mr]
\bauthor{\bsnm{El~Karoui},~\bfnm{Nicole}\binits{N.}},
  \bauthor{\bsnm{Jeanblanc},~\bfnm{Monique}\binits{M.}} \AND
  \bauthor{\bsnm{Jiao},~\bfnm{Ying}\binits{Y.}}
(\byear{2010}).
\btitle{What happens after a default: The conditional density approach}.
\bjournal{Stochastic Process. Appl.}
\bvolume{120}
\bpages{1011--1032}.
\bid{doi={10.1016/j.spa.2010.02.003}, issn={0304-4149}, mr={2639736}}
\bptnote{check year}%
\bptok{imsref}%
\end{barticle}
\endbibitem

\bibitem{ejj2}
\begin{bmisc}[mr]
\bauthor{\bsnm{El~Karoui},~\bfnm{N.}\binits{N.}},
  \bauthor{\bsnm{Jeanblanc},~\bfnm{M.}\binits{M.}},
  \bauthor{\bsnm{Jiao},~\bfnm{Y.}\binits{Y.}} \AND
  \bauthor{\bsnm{Zargari},~\bfnm{B.}\binits{B.}}
(\byear{2010}).
\bhowpublished{Conditional default probability and density. Preprint.
Available at
\url{http://people.math.jussieu.fr/\textasciitilde jiao/recherche/density.pdf}.}
\bptok{imsref}%
\end{bmisc}
\endbibitem

\bibitem{HIM}
\begin{barticle}[mr]
\bauthor{\bsnm{Hu},~\bfnm{Ying}\binits{Y.}},
  \bauthor{\bsnm{Imkeller},~\bfnm{Peter}\binits{P.}} \AND
  \bauthor{\bsnm{M{\"u}ller},~\bfnm{Matthias}\binits{M.}}
(\byear{2005}).
\btitle{Utility maximization in incomplete markets}.
\bjournal{Ann. Appl. Probab.}
\bvolume{15}
\bpages{1691--1712}.
\bid{doi={10.1214/105051605000000188}, issn={1050-5164}, mr={2152241}}
\bptok{imsref}%
\end{barticle}
\endbibitem

\bibitem{jeamatngo10}
\begin{bmisc}[auto:STB|2012/04/25|11:03:30]
\bauthor{\bsnm{Jeanblanc},~\bfnm{M.}\binits{M.}},
  \bauthor{\bsnm{Matoussi},~\bfnm{A.}\binits{A.}} \AND
  \bauthor{\bsnm{Ngoupeyou},~\bfnm{A.}\binits{A.}}
(\byear{2010}).
\bhowpublished{Quadratic backward SDE's with jumps and utility maximization of
  portfolio credit derivative. Univ. Paris Diderot}.
\bptok{imsref}%
\end{bmisc}
\endbibitem

\bibitem{Jiao09}
\begin{bmisc}[auto:STB|2012/04/25|11:03:30]
\bauthor{\bsnm{Jiao},~\bfnm{Y.}\binits{Y.}}
(\byear{2009}).
\bhowpublished{Random measure and multiple defaults. Preprint.}
\bptok{imsref}%
\end{bmisc}
\endbibitem

\bibitem{JP09}
\begin{barticle}[auto:STB|2012/04/25|11:03:30]
\bauthor{\bsnm{Jiao},~\bfnm{Y.}\binits{Y.}} \AND
  \bauthor{\bsnm{Pham},~\bfnm{H.}\binits{H.}}
(\byear{2011}).
\btitle{Optimal investment with counterparty risk: A default-density
  approach}.
\bjournal{Finance Stoch.}
\bvolume{15}
\bpages{725--753}.
\bptok{imsref}%
\end{barticle}
\endbibitem

\bibitem{kaz}
\begin{bbook}[auto:STB|2012/04/25|11:03:30]
\bauthor{\bsnm{Kazamaki},~\bfnm{M.}\binits{M.}}
(\byear{2000}).
\btitle{Continuous Exponential Martingales and BMO}.
\bseries{Lectures Notes in Math.}
\bvolume{1579}.
\bpublisher{Springer}, \baddress{Berlin}.
\bptok{imsref}%
\end{bbook}
\endbibitem

\bibitem{khalim11}
\begin{bmisc}[auto:STB|2012/04/25|11:03:30]
\bauthor{\bsnm{Kharroubi},~\bfnm{I.}\binits{I.}} \AND
  \bauthor{\bsnm{Lim},~\bfnm{T.}\binits{T.}}
(\byear{2011}).
\bhowpublished{Progressive enlargement of filtrations and backward SDEs with
  jumps. Preprint, LPMA. Available at
  \url{http://hal.archives-ouvertes.fr/ccsd-00555787/en/}.}
\bptok{imsref}%
\end{bmisc}
\endbibitem

\bibitem{kob00}
\begin{barticle}[mr]
\bauthor{\bsnm{Kobylanski},~\bfnm{Magdalena}\binits{M.}}
(\byear{2000}).
\btitle{Backward stochastic differential equations and partial differential
  equations with quadratic growth}.
\bjournal{Ann. Probab.}
\bvolume{28}
\bpages{558--602}.
\bid{doi={10.1214/aop/1019160253}, issn={0091-1798}, mr={1782267}}
\bptok{imsref}%
\end{barticle}
\endbibitem

\bibitem{LQ}
\begin{barticle}[mr]
\bauthor{\bsnm{Lim},~\bfnm{Thomas}\binits{T.}} \AND
  \bauthor{\bsnm{Quenez},~\bfnm{Marie-Claire}\binits{M.-C.}}
(\byear{2011}).
\btitle{Exponential utility maximization in an incomplete market with
  defaults}.
\bjournal{Electron. J. Probab.}
\bvolume{16}
\bpages{1434--1464}.
\bid{doi={10.1214/EJP.v16-918}, issn={1083-6489}, mr={2827466}}
\bptok{imsref}%
\end{barticle}
\endbibitem

\bibitem{Pham}
\begin{barticle}[mr]
\bauthor{\bsnm{Pham},~\bfnm{Huy{\^e}n}\binits{H.}}
(\byear{2010}).
\btitle{Stochastic control under progressive enlargement of filtrations and
  applications to multiple defaults risk management}.
\bjournal{Stochastic Process. Appl.}
\bvolume{120}
\bpages{1795--1820}.
\bid{doi={10.1016/j.spa.2010.05.003}, issn={0304-4149}, mr={2673975}}
\bptok{imsref}%
\end{barticle}
\endbibitem

\bibitem{REK}
\begin{barticle}[mr]
\bauthor{\bsnm{Rouge},~\bfnm{Richard}\binits{R.}} \AND
  \bauthor{\bsnm{El~Karoui},~\bfnm{Nicole}\binits{N.}}
(\byear{2000}).
\btitle{Pricing via utility maximization and entropy}.
\bjournal{Math. Finance}
\bvolume{10}
\bpages{259--276}.
\bid{doi={10.1111/1467-9965.00093}, issn={0960-1627}, mr={1802922}}
\bptok{imsref}%
\end{barticle}
\endbibitem

\bibitem{wag80}
\begin{bincollection}[mr]
\bauthor{\bsnm{Wagner},~\bfnm{Daniel~H.}\binits{D.~H.}}
(\byear{1980}).
\btitle{Survey of measurable selection theorems: An update}.
In \bbooktitle{Measure Theory, {O}berwolfach 1979 ({P}roc. {C}onf.,
  {O}berwolfach, 1979)}.
\bseries{Lecture Notes in Math.}
\bvolume{794}
\bpages{176--219}.
\bpublisher{Springer}, \baddress{Berlin}.
\bid{mr={0577971}}
\bptok{imsref}%
\end{bincollection}
\endbibitem

\end{thebibliography}
\end{document}